\documentclass[11pt,leqno]{amsart}
\usepackage{amsmath,amssymb,array,latexsym}
\newtheorem{thm}{Theorem}[section]
\newtheorem{lem}[thm]{Lemma}
\newtheorem{prop}[thm]{Proposition}
\newtheorem{cor}[thm]{Corollary}
\newtheorem{example}[thm]{Example}
\newtheorem{prob}[thm]{Problem}
\theoremstyle{definition}
\newtheorem*{defin}{Definition}
\newtheorem*{remark}{Remark}
\newtheorem*{claim}{Claim}
\def\hangbox to #1 #2{\vskip1pt\hangindent #1\noindent \hbox to #1{#2}$\!\!$}

\def\mod{\mathop{\rm mod}\nolimits}
\def\rW{\text{\rm W}}
\newcommand{\n}{\noindent}
\newcommand{\ds}{\displaystyle}
\newcommand{\vp}{\varepsilon}
\newcommand{\tC}{\tilde C}
\newcommand{\tX}{\tilde X}
\newcommand{\coo}{\mbox{c}_{00}}

\newcommand{\dist}{\operatorname{dist}}
\newcommand{\card}{\operatorname{card}}
\newcommand{\cof}{\operatorname{cof}}

\newcommand{\Id}{\operatorname{Id}}
\newcommand{\cT}{\mathcal{T}}
\newcommand{\cI}{\mathcal{I}}
\newcommand{\cF}{\mathcal{F}}
\newcommand{\cC}{\mathcal{C}}
\newcommand{\cD}{\mathcal{D}}
\newcommand{\cA}{\mathcal{A}}
\newcommand{\cB}{\mathcal{B}}
\newcommand{\cN}{\mathcal{N}}
\newcommand{\cS}{\mathcal{S}}
\newcommand{\bcAvp}{\overline{\mathcal{A_\vp}}}

\newcommand{\N}{\mathbb{N}}

\newcommand{\fN}{[\N]^{<\omega}}
\newcommand{\iN}{[\N]^{\omega}}

\newcommand{\iNN}{[\overline{N}]^{\omega}}
\newcommand{\kN}{[\N]^{\le k}}

\newcommand{\bN}{\overline{N}}
\newcommand{\bM}{\overline{M}}

\author{E.~Odell$^{(1)}$  and Th.~Schlumprecht$^{(2)}$}
\title{Trees and Branches in Banach Spaces}
\thanks{$^{(1),(2)}$ Research supported by NSF}
\begin{document}
\begin{abstract}
An infinite dimensional notion of asymptotic structure is considered. 
This notion is developed in terms of trees and branches on Banach spaces. 
Every countably infinite countably branching tree $\cT$ of a certain type on 
a space $X$ is presumed to have a branch with some property. 
It is shown that then $X$ can be embedded into a space with an FDD 
$(E_i)$ so that all normalized sequences in $X$ which are almost a skipped 
blocking of $(E_i)$ have that property. 
As an application of our work we prove that if $X$ is a separable reflexive 
Banach space and for some $1<p<\infty$ and $C<\infty$ every weakly null 
tree $\cT$ on the sphere of 
$X$ has a branch $C$-equivalent to the unit vector basis of 
$\ell_p$, then for all $\vp>0$, 
there exists a finite codimensional subspace of $X$ which 
$C^2+\vp$ embeds into the $\ell_p$ sum 
of finite dimensional spaces. 
\end{abstract}
\maketitle
\markboth{E. ODELL AND TH. SCHLUMPRECHT}{TREES AND BRANCHES IN BANACH SPACES}

\baselineskip=15pt

\section{Introduction}\label{S0}

A recurrent theme in Banach space theory takes the following form. 
One has some property $(P)$ and one assumes that in a given separable infinite 
dimensional Banach space $X$, every normalized weakly null sequence (or 
perhaps every normalized block basis of a given basis for $X$) admits a 
subsequence with $(P)$. 
One then tries to deduce that $X$ has some other property $(Q)$. 
In this paper we consider a stronger hypothesis on $X$. 
Namely that every countably infinitely branching tree of 
$\omega$-levels of some  type (e.g., the successors of every node are a 
normalized weakly null sequence or perhaps a block basis of some FDD) 
admits a branch with $(P)$. 
As we show this is sometimes the proper hypothesis to conclude that $X$ 
has $(Q)$. 

An example of this type is given in Theorem~4.1 where the following is 
proved: 
If $X$ is reflexive and there exists $1<p<\infty$ and $C<\infty$ so that 
every normalized weakly null tree in $X$ admits a branch $C$-equivalent to 
the unit vector basis of $\ell_p$ then for all $\vp>0$ 
there exists a finite codimensional subspace of $X$ which 
$C^2+\vp$-embeds into some space $(\sum F_i)_p$, an $\ell_p$-sum of 
finite dimensional spaces. 
Hence this characterizes when a reflexive space embeds into such a sum. 

The motivation for working with branches of trees in place of subsequences 
comes from the notion of asymptotic structure (\cite{MT}, \cite{MMT}),  
the recent paper of N.J.~Kalton \cite{K} and \cite{KOS}. 
In its simplest version suppose $X$ has an FDD $(E_i)$ and let $k\in\N$. 
Then the $k^{th}$-asymptotic space of $X$ with 
respect to $(E_i)$ may be described
as the smallest closed set $C_k$ of normalized bases of length $k$ with 
the property that every countably infinitely branching tree of $k$ levels 
in $S_X$ whose nodes are all block bases of $(E_n)$ must admit for every 
$\vp>0$ a branch $1+\vp$-equivalent to some member of $C_k$. 

Moreover given $\vp_n\downarrow0$ one can then block $(E_n)$ into an FDD 
$(F_n)$ with the property that for all $k$ any normalized skipped block 
basis $(x_i)_1^k$ of $(F_n)_{n=k}^\infty$ is $1+\vp_m$-equivalent to 
a member of $C_k$ \cite{KOS}. 
We cannot achieve this in the infinite setting, $k=\omega$. 
There is in general no unique infinite asymptotic structure, $C_\omega$. 
However if $C$ is big enough so that every such $\omega$-level tree has a 
branch in $C$ then one can produce for $\vp>0$ a blocking $(F_n)$ of $(E_n)$ 
so that all normalized skipped block bases of $(F_n)$ starting after $F_1$ 
are in $\overline{C_{\vp}}$, the pointwise closure (in the product topology 
of the discrete topology on $S_X$) of 
$\frac{\vp}{2^n}$-perturbations of elements of $C$. 
This is done in section~3.  
(We note that an in between 
ordinal notion of asymptotic structure for $\alpha<\omega_1$ has been 
considered in \cite{W}, using the generalized Schreier sets $S_\alpha$.) 

Actually we need to study more general forms of asymptotic structure than 
that w.r.t. an FDD. 
We consider the version where one uses arbitrary finite codimensional 
subspaces rather than just the tail subspaces of a given FDD. 
While this version is coordinate free we show in section~3 that one may 
embed $X$ into a space with an FDD in such a way that the two notions coincide. 
Section~2 contains our preliminary work and terminology. 
In section~5 we apply our results to the more general notion of V.D.~Milman's 
\cite{Mi}  spectra of a function. 
We are indebted to W.B.~Johnson for showing us the proof of Lemma~3.1.

\section{Games in a Banach space $X$}\label{S1}

Assume that $X$ is a separable Banach space of infinite dimension. The set of all subspaces of $X$
having  finite codimension is denoted by $\cof(X)$.  
 $S_X^{\,\omega}$ and $ S_X^k$, $k\in \N$, denote the
set of all infinite sequences in $S_X$, the unit sphere of $X$, respectively 
all sequences in $S_X$ of length $k$.

For a  set $\cA\subset S_X^{\,\omega}$ or  $\cA\subset S_X^k$ 
we consider the following
{\em $\cA$-game\/} between two players, having infinitely many,  
respectively $k$, rounds: 
\begin{align*}
&\text{Player I chooses } Y_1\in\cof(X)\\
&\text{Player II chooses }y_1\in S_{Y_1}\\
&\text{Player I chooses }Y_2\in\cof(X)\\
&\text{Player II chooses }y_2\in S_{Y_2}\\
&\ldots
\end{align*}
Player I wins if the resulting sequence $(y_i)$ is in $\cA$.
  
Note that by replacing a set $\cA\subset S^k_X$, $k\in \N$, 
by $\cA\times S_X^{\,\omega}$, we need 
only consider games with infinitely many steps.

We say that {\em Player I has a winning strategy in the $\cA$-game\/} 
if the following condition, $\text{W}_I(\cA)$ holds.
\begin{equation}\tag*{$(\rW_I(\cA))$}
\left\{ \begin{array}{cl}
&\text{There is a family of finite codimensional subspaces of $X$}\\
\noalign{\vskip4pt}
&\qquad \bigl( Y_{(x_1,x_2,\ldots x_\ell)}\bigr)_{(x_1,x_2,\ldots x_\ell)\in 
\bigcup_{j=0}^\infty S^j_X}\ , \qquad S^0_X=\{\emptyset\}\ ,\\
\noalign{\vskip4pt}
&\text{indexed over all finite sequences in $S_X$, so that:}\\
&\text{If $(x_n)_{n\in\N}$  satisfies the following recursive condition:}\\
\noalign{\vskip4pt}
&(1)\qquad x_1\in S_{Y_\emptyset}, \text{ and, for }n\ge 2,\ 
x_n\in S_{\ds Y_{(x_1,\ldots x_{n-1})}}\ ,\\
\noalign{\vskip4pt}
&\text{then $(x_n)\in \cA$.} 
\end{array} \right.
\end{equation}
\addtocounter{equation}{1}

The following Proposition can be deduced immediately from the 
definition of $(\text{W}_{I}(\cA))$. 

\begin{prop}\label{P0}
The set of all $A\subset S^{^\omega}_X$ for which Player~I has 
a winning strategy is closed with respect to taking finite intersections.
\end{prop}

\n Similarly, we say that Player II has a winning strategy if
 
\begin{equation}\tag*{$(\text{W}_{II}(\cA))$}
\left\{ \begin{array}{cl}
&\text{There is a family in $S_X$}\\
\noalign{\vskip4pt}
&\qquad\qquad\qquad
\bigl(x_{(Y_1,Y_2,\ldots Y_\ell)}\bigr)_{(Y_1,Y_2,\ldots Y_\ell)
\in\bigcup_{j=1}^\infty \cof^{\,j}(X)}\ ,\\
\noalign{\vskip4pt}
&\text{indexed over all finite sequences in $\cof(X)$ 
(of length at least 1) so that}\\
\noalign{\vskip4pt}
&(2)\qquad  x_{(Y_1,Y_2,\ldots Y_\ell)}\in S_{Y_\ell}
\text{ if }\ell\in\N\text{ and }Y_1,\ldots ,Y_\ell\in\cof(X),\ \text{ and}\\
\noalign{\vskip4pt}
&(3)\qquad \text{for every sequence }(Y_i)_{i\in\N}\subset\cof(X), 
\bigl(x_{(Y_1,Y_2,\ldots Y_i)}\bigr)_{i=1}^\infty\not\in\cA.    
\end{array} \right.
\end{equation}
\addtocounter{equation}{2}

\begin{remark} 
Informally $(\text{W}_{I}(\cA))$ means the following:
$$\exists Y_1\in\cof(X)\forall y_1\in S_{Y_1}\exists Y_2\in\cof(X)
\forall y_2\in S_{Y_2}\ldots\text{ so that }
(y_i)\in\cA.$$
Since this is an infinite phrase (unless we considered a game 
of finitely many draws), it has to be defined
in a more formal way as it was done in  $(\text{W}_{I}(\cA))$. 
 
It is not true in general that an $\cA$-game is determined, i.e.,
that either Player~I or Player~II has a winning strategy. 
Note that this would mean that if the  above infinite phrase is false then 
we can formally negate it. 
  
{From} a result of D.~A.~Martin [Ma] it follows that if $\cA$ is a 
Borel set with respect to the product topology
of the discrete topology in $S_X$ then the $\cA$-game is determined. 
We actually will only need a special case of this theorem which is much 
easier (see [GS] or section 1 of [Ma]). 
\end{remark}
  
\begin{prop}\label{P1}
For every $\cA\subset S_X^{\,\omega}$  
$(\text{W}_{I}(\cA))$ and  $(\text{W}_{II}(\cA))$ are mutually exclusive
and if $\cA$ is closed with respect to the product  of the discrete topology,
then it follows that the failure of  $(\text{W}_{I}(\cA))$ 
implies $(\text{W}_{II}(\cA))$.
  
We furthermore note that both statements remain 
true if we change the game to a game
in which Player I has to choose his spaces  among some given subset 
$\Gamma\subset\cof(X)$ and/or
Player II has to choose his vectors among a subset 
$D\subset S_X$ or can choose
his vector in some neighborhood of $S_{Y_n}$, 
with $Y_n$ being the $n$-th choice of Player I.
\end{prop}

For a more detailed description of these variations of the 
$\cA$-game we refer to Proposition
\ref{P2}, where we discuss the existence of winning strategies.
In that Proposition we will show that we can reduce the 
game into a game in which Player I,
assuming he has a winning strategy, can  determine a countable collection 
of finite codimensional spaces before the game starts,
then make his choices among this countable collection and still win the game.
  
We need the following notion of trees and some terminology.
  
\begin{defin} 
$[\N]^{<\omega}$ denotes the set of nonempty finite subsets of $\N$ and 
$[\N]^{\le k}$ denotes the nonempty subsets of $\N$ of cardinality at 
most $k$. 
These are regarded as countably branching trees of infinite length, 
respectively, of length $k$, under the order $A\le B$ if $A$ is an initial 
segment of $B$. 
A {\em countably branching tree of infinite length in\/} $S_X$ is a family
$(x_A)_{A\in\fN}$ in $S_X$, where the order is that induced by $\fN$.  
\end{defin}
  
Similarly a 
{\em countably branching tree of  length $k\in\N$ in\/} $S_X$ is a family
$(x_A)_{A\in\kN}$ in $S_X$. 
  
Since these are the only kinds of trees we will consider we will simply 
refer to them as 
{\em trees of infinite or finite length in\/} $S_X$.
 
If  $(x_A)_{A\in\fN}$ or $(x_A)_{A\in\kN}$ is a tree and 
$A\in\fN\cup\{\emptyset\}$, or 
$A\in[\N]^{\le k-1}\cup\{\emptyset\}$ respectively,
we call  the sequence $(x_{A\cup\{n\}})_{n>\max A}$
{\em the $A$-node\/} of that tree.
 
If $(n_i)$ is an increasing sequence in $\N$ of infinite length, 
respectively of length $k$, we call the sequence
$(x_{\{n_1,\ldots n_i\}})_{i=1}^\infty$, respectively  
$(x_{\{n_1,\ldots n_i\}})_{i=1}^k$, a {\em branch of the tree\/}.
 
Assume that $(x_A)_{A\in\fN}$ or $(x_A)_{A\in\kN}$ 
is a tree of infinite length or length $k$, respectively, 
and $\cI\subset\fN$, or $\cI\subset\kN$ has the following property:
\begin{enumerate}
\item[a)] $\cI$ is hereditary, i.e.,
if $A\in\cI$, and $\emptyset\ne B$ is an initial 
segment of $A$ then $B\in\cI$.
\item[b)] Assume that $A\in\cI\cup\{\emptyset\}$,
and that $\card(A)<k$, if we consider the case of a tree of length $k$. Then
there are infinitely many direct
successors of $A$ in $\cI$, i.e., 
the set $\{n\in\N:A\cup\{n\}\in\cI\}$ is infinite.
\end{enumerate}
Then we call the family $(x_A)_{A\in\cI}$ a {\em subtree of\/} $(x_A)$.
Note that in that case we can relabel the family $(x_A)_{A\in\cI}$ as 
a tree $(y_A)_{A\in\fN}$  or $(y_A)_{A\in\kN}$, respectively, so that every 
node and every branch of $(x_A)_{A\in\cI}$ is node or branch, 
respectively, of $(y_A)$ and vice versa.
   
If $(Y_n)$ is a decreasing sequence of finite codimensional subspaces of $X$, 
we call a tree $(x_A)$
(indexed over $\fN$ or $\kN$) a $(Y_n)$-{\em block-tree\/} if for every
$A\in\fN$, respectively every $A\in\kN$, $x_A\in S_{Y_{\max A}}$.
   
Let $\delta_i\in(0,1]$, for $i\in\N$, $\delta_i\searrow 0$. 
We call a   tree $(x_A)_{A\in\fN}$ of infinite length in $S_X$ a 
{\em $(\delta_i)$- approximation of a $(Y_n)$-block tree\/}, if 
$$\dist(x_A, S_{Y_{\max A}})<\delta_{\card A}, \text{ whenever }A\in\fN$$
  
If $\cT$ is a topology on $X$ (for example the weak topology), we call a 
tree $\cT${-null} if every node is a $\cT$-null sequence. 
 
\begin{remark} 
For a sequence $(x_n)\subset X$  we can define a 
tree $(x_A)_{A\in\fN}$, by setting $x_A:=x_{\max A}$, for $A\in\fN$. 
Note that then the set of all subsequences of $(x_n)$ coincides
with  the set of all branches of $(x_A)_{A\in\fN}$.
\end{remark}
  
We will be interested in conditions of the following form 
and relate them to the existence of
winning strategies of the above discussed games.
  
\begin{enumerate}
\item[]
{\em Assume that all trees all of whose nodes have a certain property (A) 
(for example being weakly null), have a branch
with a certain property (B) (for example being equivalent to 
the unit vector basis of $\ell_p$).}
\end{enumerate}
{From} the above,  such a  condition is a strengthening of the 
following assumption:
\begin{enumerate}
\item[]
{\em All normalized  sequences  having property (A) have a 
subsequence with property (B).} 
\end{enumerate}

Continuing with our notation, 
if $\cA\subset S^{\,\omega}_X$ and $\vp>0$, we let 
$$\cA_\vp=\bigl\{(x_i)\subset S_X:\exists (y_i)\in\cA, \|x_i-y_i\|<\vp/2^i 
\text{ for all } i\in\N\bigr\}$$
and let $\bcAvp$ be the closure of $\cA_\vp$  with respect to the product 
of the discrete topology.
We note that for $\vp,\delta >0$ 
\begin{equation}\label{E4}
\overline{(\bcAvp)_\delta}\subset\overline{\cA_{\vp+\delta}}\ .
\end{equation}
  
If $Y\in \cof (X)$ and $\delta >0$ then 
$$(S_Y)_\delta = \{x\in S_X :\|x-y\| <\delta \text{ for some } 
y\in S_Y\}\ .$$

Let $\vp >0$, $\Gamma \subseteq \cof (X)$ and $D\subseteq S_X$. 
We define what it means to say Player~I has a winning strategy 
for $\cA\subset S_X^\omega$ given that Player~I 
can only choose $Y\in\Gamma$ or that II can only choose elements of $D$.
\begin{equation}\tag*{$(\text{W}_I(\cA,\Gamma,\vp))$} 
\left\{ \begin{array}{cl} 
&\text{There exists a family}\\
\noalign{\vskip4pt}
&\qquad\qquad\qquad 
\bigl(Y_{(x_1,x_2,\ldots x_\ell)}\bigr)_{(x_1,x_2,\ldots x_\ell)
\in\bigcup_{j=0}^\infty S^j_X}\subset\Gamma\ ,\\
\noalign{\vskip4pt}
&\text{so that for every sequence $(x_n)_{n\in\N}$ satisfying the following}\\ 
&\text{recursive condition:}\\
\noalign{\vskip4pt}
&(5)\qquad\quad x_1\in (S_{Y_\emptyset})_{\vp/2}, \text{ and, for }n\ge 2,
\ x_n\in (S_{Y_{(x_1,\ldots x_{n-1})}})_{\vp/2^n}\\
\noalign{\vskip4pt}
&\text{one has $(x_n)\in\cA$.}
\end{array} \right.
\end{equation}
\addtocounter{equation}{1}
          
\begin{remark} 
It is easy to see by \eqref{E4} that for any $\vp,\delta>0$, 
$$(\text{W}_I(\overline{\cA_{\vp}},\{Y_n\},\vp))\Rightarrow  
(\text{W}_I(\overline{\cA_{\vp+\delta}},\{\tilde Y_n\},\vp))$$
whenever $\{\tilde Y_n\}\subseteq \cof (X)$ is a {\em refinement of\/} 
$\{Y_n\}$, by which we mean that 
$$\forall\ Y\in\{Y_n\}\ \forall\ \delta>0\ \exists\ \tilde Y\in\{\tilde Y_n\}
\text{ with } S_{\tilde Y}\subset (S_Y)_\delta\ .$$   
\begin{equation}\tag*{$(\text{W}_I(\cA,D,\vp))$} 
\left\{ \begin{array}{cl}
&\text{There is a family}\\
\noalign{\vskip4pt}
&\qquad\qquad
\bigl(Y^{(\vp)}_{(x_1,\ldots x_\ell)}\bigr)_{(x_1,\ldots x_\ell)
\in\bigcup_{j=0}^\infty D^j}\subset \cof(X)\ ,\\
\noalign{\vskip4pt}
&\text{so that for any sequence $(x_n)$, such that $x_{n}\in D$, and}\\
\noalign{\vskip4pt}
&\qquad\qquad x_{n}\in ( S_{Y^{(\vp)}_{(x_1,\ldots x_{n-1})}})_{\vp/2^{n}}\ ,
\text{ $n=1,2,\ldots$,}\\
\noalign{\vskip4pt}
&\text{one has $(x_n)\in \cA$.}
\end{array} \hskip.5truein\right.
\end{equation}
\end{remark}

\begin{prop}\label{P2} 
  
\begin{enumerate}
\item If $\cB$ is a countable collection of subsets of  
$ S^{\,\omega}_X$, then there is a decreasing sequence  $(Y_n)$ 
in $\cof(X)$ so that the following are equivalent for each $\cA\in\cB$
\begin{enumerate}
\item[a)] $\forall\vp>0\quad (\text{W}_I(\,\overline{\cA_\vp})\,)$.
\item[b)] $\forall\vp>0\quad (\text{W}_I(\overline{\,\cA_{2\vp}},\{Y_n\},\vp))$.
\item[c)]  For every $\vp>0$ every $(\vp/2^n)$-approximation to a 
$(Y_n)$ block tree of infinite length in $S_X$  has a branch in 
$\overline{\cA_{\vp}}$.
\item[d)]  For every $\vp>0$ every  $(Y_n)$ block tree of infinite length in 
$S_X$  has a branch in $\overline{\cA_{2\vp}}$.
\end{enumerate}
  
\item  
If $X$ has a separable dual, then $(Y_n)\subset\cof(X)$ can be chosen so that
the equivalences in 1.\  
hold for all subsets $\cA\subset S^{\,\omega}_X$.
In that case it follows that
for any $\cA\subset S_X^{\,\omega}$ that (1)(a) is equivalent to       
\begin{enumerate}
\item[e)] For every $\vp>0$ every weakly null tree of infinite length in $S_X$  
has a branch in $\overline{\cA_\vp}.$
\end{enumerate}
\end{enumerate}
\end{prop}
  
\begin{proof}[Proof of Proposition \ref{P2}:] 
Let $D$ be a countable dense set in $S_X$.
Using (4) we note that for any
$\cA\subset S^{\,\omega}_X$  and  any $\vp>0$ it follows that 
\begin{equation}\label{E4a}
(\text{W}_I(\,\bcAvp\,))\Rightarrow
(\text{W}_I(\overline{\,\cA_{2\vp}},D,\vp))\ .
\end{equation}
  
Assuming now that for all $\vp>0$ the condition 
$(\text{W}_I(\,\overline{\cA_{2\vp}},D,\vp))$ 
is satisfied we can choose a countable subset of $\cof (X)$,  
\begin{equation}\label{E4aa}
\Gamma_\cA=\big\{ Y^{(\vp)}_{(x_1,\ldots x_\ell)} : 
\vp>0\text{ rational },x_n\in D\text{ and } x_n\in 
(Y_{(x_1,\ldots,x_{n-1})}^{(\vp)})_{\vp/2^n}
\text{ for } n\in\N\bigr\}\ , 
\end{equation}
and observe that 
\begin{align}\label{E4b}  
\forall\vp>0\quad (\text{W}_I(\,\overline{\cA_{2\vp}},D,\vp)) \implies 
 &\text{there exists a countable }
 \Gamma\subset\cof(X) \text{ so that }\\ 
 &\forall\vp>0,
\quad (\text{W}_I(\,\overline{\cA_{2\vp}},\Gamma,D,\vp)).\notag
\end{align} 
where $(\text{W}_I(\,\overline{\cA_{2\vp}},\Gamma,D,\vp))$ is defined just 
like   $(\text{W}_I(\,\overline{\cA_{2\vp}},\Gamma,\vp))$ with the
difference that the family $\bigl(Y_{(x_1,x_2,\ldots x_\ell)}\bigr)$ is indexed
over $\bigcup_{j=0}^\infty D^j$
   
Using standard approximation arguments and the fact that 
$D$ is dense in $S_X$ we observe for
any $\Gamma\subset\cof(X)$ and any  $\cA\subset S^{\,\omega}$
\begin{equation}\label{E4c}
(\text{W}_I(\,\overline{\cA_{2\vp}},\Gamma,D,\vp))
\Rightarrow (\text{W}_I(\,\overline{\cA_{3\vp}},\Gamma,\vp))
\Rightarrow(\text{W}_I(\,\overline{\cA_{3\vp}}\,))\ .
\end{equation}
    
Finally assume that $\tilde\Gamma \subset \cof(X)$ is a  refinement 
of $\Gamma\subset\cof(X)$.
Then by \eqref{E4} it follows for $\vp>0$ that
\begin{equation}\label{E4d}
(\text{W}_I(\,\bcAvp,\Gamma))\Rightarrow (\text{W}_I
(\,\overline{\cA_{2\vp}},\tilde\Gamma)).
\end{equation}
     
Let $\cB$ be any countable collection of subsets of $S_X^{\,\omega}$. 
For $\cA\in\cB$, if for all $\vp>0$ (W$_I(\,\bcAvp\,)$) is true let $\Gamma_A$ 
be as in \eqref{E4aa}, and, otherwise, we set $\Gamma_{\cA}=\{X\}$. 
Since $\bigcup_{\cA\in\cB}\Gamma_{\cA}$ is countable we can choose
a decreasing sequence $(Y_n)\subset\cof$ which is a refinement of 
$\bigcup_{\cA\in\cB}\Gamma_{\cA}$.
      
\n {From} \eqref{E4a}--\eqref{E4d}
we deduce that for all $\cA\in\cB$
$$\forall \vp>0\quad \text{W}_I(\,\bcAvp\,)\iff\forall\vp>0\quad 
\text{W}_I(\,\overline{\cA_{2\vp}},\{Y_n\},\vp)\ .$$
Now $\text{W}_I(\,\overline{\cA_{2\vp}},\{Y_n\},\vp)$ says 
that Player I has in the 
$\overline{\cA_{2\vp}} $-game a winning strategy,
even if he has to choose his finite codimensional subspaces among $\{Y_n\}$,
and even if Player II ``can cheat a little bit'' by choosing his vectors 
in $(S_{Y_n})_{\vp/2^n}$.
{From} Proposition~\ref{P1} we deduce that this is equivalent 
to the condition that 
Player II does not have a winning strategy which means that every 
$(\vp/2^n)$ approximation to a
$(Y_n)$-block-tree has a branch in $\overline{\cA_{2\vp}}$.
       
We therefore have proven the equivalence of (a), (b) and (c). 
Note also that (c)$\Rightarrow$(d) is trivial and
since (d) means that Player~II has no winning strategy even if Player~I 
has to choose form the set $\{Y_n\}$ it follows that (d) implies (a). 
      
In order to prove the second part of the Proposition we note that in the case
that $X$ has a separable dual we can find a {\em universal\/} 
countable refinement, i.e., a countable refinement of the whole set $\cof(X)$. 
Indeed, choose a dense sequence $(\xi^*_n)$ in $S_{X^*}$ and let 
$$Y_n=\cN(\xi^*_1, \xi^*_2,\ldots , \ldots\xi^*_n)
=\{x\in X:\forall i\in\{1,\ldots n\}\quad\xi^*_i(x)=0\}\ .$$
       
Secondly note that in this case every $(Y_n)$-block-tree is weakly null, 
and, conversely, that for $\delta_i\searrow 0$, every weakly null tree 
$(x_A)_{A\in\fN}$ has a subtree $(y_A)_{A\in\fN}$ which
is a  {\em$(\delta_i)$-approximation\/} of a  
$(Y_n)$-block-tree.
\end{proof}
      
\section{A fundamental combinatorical result}\label{S2}
 
For the games in $X$, introduced in Section \ref{S1}, we want 
to discuss how a winning strategy of Player I or Player II can be
formulated in terms of a coordinate system on $X$.
 
Recall that a Banach space $Z$ has an FDD $(F_i)$, where, 
for $i\in\N$, $F_i$ is a finite dimensional
subspace of $Z$, if every $z\in Z$  can be written in a unique way as
$z=\sum_{i=1}^\infty z_i$ with $z_i\in F_i$, for all $i\in\N$.
In this case we write $Z=\oplus_{i=1}^\infty F_i$ and denote by 
$\coo(\oplus_{i=1}^\infty F_i)$
the dense linear subspace of $Z$ consisting of all 
finite linear combinations of vectors $x_i$, $x_i\in F_i$.
For $m\le n$   we denote by $P_{\oplus_{i=m}^n F_i}$  
the canonical projection form $Z$ onto $\oplus_{i=m}^n F_i$.
    
Using a result of  W.~B.~Johnson, H.~Rosenthal and M.~ Zippin [JRZ]
we derive the following Lemma.
 
\begin{lem}\label{L1}
Let $(Y_n)$ be a decreasing sequence of subspaces of $X$, each  
having finite codimension.
Then $X$ is isometrically embeddable into a space $Z$
having an FDD $(E_i)$ so that
(we identify $X$ with its isometric image in $Z$)
\begin{enumerate}
\item[a)]   $\coo(\oplus_{i=1}^\infty E_i)\cap X$ is dense in $X$. 
\item[b)]  For every $n\in\N$ the finite codimensional 
subspace $X_n=\oplus_{i=n+1}^\infty E_i\cap X $ 
is contained in $Y_n$.
\item[c)]  There is a $c>0$, so that for every $n\in\N$, there is a finite
set $D_n\subset S_{\ds\oplus_{i=1}^n E_i^*}$ such that 
whenever $x\in X$ 
\begin{equation}\label{E4e}
\|x\|_{X/Y_n}=\inf_{y\in Y_n}\|x-y\|
\le c\max_{w^*\in D_n} w^*(x) \ .
\end{equation}
\end{enumerate}
{From} (a) it follows that   $\coo (\oplus_{i=n+1}^\infty E_i)\cap X$ 
is a dense linear subspace of $X_n$.

Moreover if $X$ has a separable dual $(E_i)$ can be chosen to be shrinking
(every normalized block sequence in $Z$ with respect to $(E_i)$ converges
weakly to $0$, or, equivalently, $Z^* =\oplus_{i=1}^\infty E_i^*$), 
and if $X$ is reflexive $Z$ can also be chosen to be reflexive.  
\end{lem}
   
\begin{remark} 
We will  prove that $X$ is isomorphic to a space 
$\tilde X$ having above properties.
Then we consider on $\tilde X$ the norm, $\| I(\cdot)\|_X$, where 
$I: \tilde X\to X$ is an isomorphism,
and extend this norm to all of $Z$. 
We might loose monotonicity, or bimonotonicity, and we will not be able to 
assume that the constant $c$ in (c) can be chosen close to the value~1. 
But for later purposes we are more interested in an isometric embedding.
\end{remark}
   
\begin{proof}[Proof of Lemma \ref{L1}]
We consider the following three cases.
If $X$ is a reflexive space we can choose according to [Z] a reflexive space
$Z$ with an FDD  $(F_i)$ which contains $X$.
If the dual $X^*$ is separable we can use again a result in [Z]  and
choose a space $Z$ having a shrinking FDD $(F_i)$.
In the general case we choose $Z$ to be a C$(K)$-space  containing $X$,
$K$ compact and metric 
(for example $K=B_{X^*}$ endowed with the $w^*$-topology) and choose
an FDD $(F_i)$ for $Z$.
  
We first write $Y_n$ as the null space $\cN(U_n)$ of a finite dimensional
space $U_n\subset X^*$ . 
We choose a finite set in $S_{U_n}$, which norms all elements of
$X/Y_n$ up to a factor $1/2$ and  choose for each element of this
set a Hahn-Banach extension to an element  in $Z^*$. 
We denote the set of all extensions by  $D_n$ and 
let $V_n$ be the finite dimensional
subspace of $Z^*$ generated by $D_n$.  
We will produce an FDD $(E_i)$ for $Z$ so that 
$D_n\subset \oplus_{i=1}^n E_i^*$. 
Hence (c) will hold. 
     
Now   
\begin{equation}\label{E5}
Y_n=\cN(V_n)\cap X, \text{ with }V_n\subset Z^*,\text{ and }\dim(V_n)<\infty\ .
\end{equation}
Secondly we choose a subspace  $\tilde W_n\subset X$, 
$\dim(\tilde W_n)=\dim(U_n)<\infty$,
so that $X$ is a complemented sum of $Y_n$ and $\tilde W_n$, 
$X= Y_n \oplus \tilde W_n$. 
Note that in general we do not have  control over the norm 
of the projection onto $Y_n$.  
Given a dense countable subset $(\xi_n)$ in $S_X$, we inflate
$\tilde W_i$ to $W_i=\text{span}(\tilde W_i \cup \{\xi_1,\ldots\xi_i\})$.
Thus the closure of $\bigcup_{i=1}^\infty W_i$ is $X$.
    
Then we choose as follows a separable subspace $\tilde Z$ of $Z^*$
which is 1-complemented in $Z^*$,  $Z$-norming, and contains
all the spaces $V_n$, $n\in\N$.
In the case that $X$ has a separable dual (thus also $Z^*$
is separable) we simply take $\tilde Z=Z^*$. In the general
case we let $\tilde Z$ be a separable $L_1$-space containing
a $Z$-norming set,  all the spaces $V_n$, and all the spaces $F_n^*$
(considered as subspaces of $Z^*$).  
  
For $n\in\N$ let $P_n: Z\to \oplus_{i=1}^n F_i$ be the projection from $Z$ onto
$\oplus_{i=1}^n F_i$, and let $T_n:Z^*\to \tilde Z$ be 
the adjoint $P^*_n$ if $X^*$ is separable. 
In the general case we choose $(T_n)$ to be a sequence of
projections of norm~1 from $Z^*$ onto a finite dimensional subspace of 
$\tilde Z$ with the property $T_1(Z^*)\subset T_2(Z^*)\subset T_3(Z^*)\ldots$
so that $\bigcup_n T_n(Z^*)$ is dense in $\tilde Z$ (as a separable 
L$_1$-space $\tilde Z$ is complemented in $Z^*$ and  has an FDD).
      
We are now in the situation of Lemma 4.2 of [JRZ], i.e., the
following statements hold:
\begin{gather}
P_n^*(Z^*)\subset \tilde Z\text{ and }T_n(Z^*)\subset \tilde Z\ ,\label{E6}\\
\lim_{n\to\infty} P_n(z)=z,\quad\lim_{n\to\infty} T_n(y^*)=y^*
\text{ for all }z\in Z, y^*\in \tilde Z,\text{ and }\label{E7}\\
K:=\sup_n\|T_n\|\vee\sup_n\|P_n\|<\infty.\label{E8}
\end{gather}
We conclude from Lemma 4.2 in [JRZ] that:
\begin{enumerate}
\item[$(*)$] 
Let $E$ and $F$ be finite dimensional 
subspaces of $X$ and $\tilde Z$ respectively.
Then there is a projection $Q$ on $Z$ 
with finite dimensional range  so that  
the following three conditions \eqref{E9}, \eqref{E10} and \eqref{E11} hold
\end{enumerate}
\begin{gather}
Q|_E=\Id|_E\text{ and }Q^*|_F=\Id|_F\label{E9}\\
Q^*(Z^*)\subset \tilde Z\label{E10}\\
\|Q\|\le 4(K+K^2)\label{E11}
\end{gather}
  
Using $(*)$  we can proceed as in the proof of
Theorem~4.1 in [JRZ] to inductively define for each $n\in\N$
a finite dimensional projection
$(Q_n)$ on $Z$ so that for all $1\le i,j\le n$ 
\begin{gather}
Q_iQ_j=Q_jQ_i=Q_{i\wedge j}\ ,\label{E12}\\
Q_i(X)\supset\bigcup_{s=1}^i W_s\ ,\label{E13}\\
\tilde Z\supset Q_i^*(Z^*)\supset\bigcup_{s=1}^i V_s
\text{ (in particular $D_i\subset Q^*_i(Z^*)$), and}\label{E14}\\ 
\|Q_i\|\le 4(K+K^2)\ .\label{E15}
\end{gather}
Indeed, for $n=1$ we apply $(*)$ to $E=W_1$ and $F=V_1$. 
If $Q_1,Q_2,\ldots Q_{n-1}$ are
chosen we apply $(*)$ to $E=[Q_{n-1}(Z)\cup W_n]$ and 
$F=\text{span}(Q_{n-1}^*(Z^*)\cup V_n)$.
We deduce \eqref{E13}, \eqref{E14} and \eqref{E15}, and we observe that
for $i<n$, $Q_n\circ Q_i=Q_i$ and $Q_n^*\circ Q_i^*=Q_i^*$.
Since for $z\in Z$ and $z^*\in Z^*$ the second equality implies that
$$\langle Q_i\circ Q_n(z),z^*\rangle
=\langle z,Q_n^*Q_i^*(z^*)\rangle 
=\langle z,Q_n^*(z^*)\rangle 
=\langle Q_i(z),z^*\rangle\ ,$$
we also deduce that $Q_i\circ Q_n=Q_i$.

Now we let $E_i=(Q_i-Q_{i-1})(Z)$ ($Q_0=0$) and deduce from 
\eqref{E12} and \eqref{E15}, that $(E_i)$ is an
FDD of a subspace of $Z$ which, by \eqref{E13} still contains $X$. 
\eqref{E13} also implies that $\coo(\oplus F_i)\cap X$ is dense in $X$. 
Putting $X_n=\oplus_{i=n+1}^\infty F_i\cap X$,
we note that for $x\in X_n$ and $z^*\in V_n$ it follows from \eqref{E14} that 
$\langle z^*,x\rangle = \langle Q_n^*(z^*),x\rangle 
=\langle z^*,Q_n(x)\rangle=0$, and thus, that $X_n\subset\cN(V_n)\cap X=Y_n$.
 
We also deduce that for $n\in\N$, $\coo(\oplus_{i=n+1}^\infty F_i)\cap X$ 
is dense in $X_n$ using 
the following Lemma which seems to be folklore.
\end{proof}
 
\begin{lem}\label{L2}
If $Y$ is a linear and dense subspace of $X$ 
and $\tilde X$ has finite codimension in $X$, 
then $\tilde X\cap Y$ is also dense in $\tilde X$.
\end{lem}

\begin{proof}
Let $F\subset X$ be a subspace of dimension $\dim(X/\tilde X)$, admitting 
a continuous  projection $Q:X\to F$, so that $(\Id-Q)(X)=\tilde X$.

Let $x\in\tilde X$. By assumption we find a sequence $(y_n)\subset Y$ 
converging to $x$.
Let $V$ be the (finite dimensional) vector space generated by 
$(Q(y_n))_{n\in\N}$ and choose a basis of $V$ of the form 
$\{Q(y_{n_1}),\ldots Q(y_{n_\ell})\}$. 
We represent each vector $Q(y_n)$ as
$$ Q(y_n)=\sum_{i=1}^\ell \lambda_i^{(n)}Q(y_{n_i})\ ,$$
and put $x_n=y_n- \sum_{i=1}^\ell  \lambda_i^{(n)}y_{n_i}$.
Note that $x_n\in Y$ and that $Q(x_n)=0$, for all $n\in\N$.
Furthermore it follows that since $\lim_{n\to\infty}\|Q(y_n)\|=0$ 
and since $(Q(y_{n_i}))_{i=1}^\ell$ is
basis of $V$, that $\lim_{n\to\infty} \lambda^{(n)}_i=0$ 
for all $1\le i\le\ell$.
Therefore it follows that 
$\lim_{n\to\infty} x_n=\lim_{n\to\infty} y_n=x$.
\end{proof}
  
We are now ready to state and to prove the main result of this section. 
If a Banach space  $Z$ has an FDD $(E_i)$, we will call a sequence $(z_i)$ 
in $Z$ {\em a block sequence with respect to\/} $(E_i)$, 
if for some $0=k_0<k_1<k_2\ldots$ 
for every $i\in\N$, $z_i\in\oplus_{j=1+k_{i-1}}^{k_i} E_j$. 
We will call a tree $(z_A)_{A\in\fN}$ or $(z_A)_{A\in\kN}$ in $S_Z$ a 
{\em $(E_i)$-block tree\/}
if every node is a block sequence with respect to  $(E_i)$.
In a similar way given $\delta_n\downarrow 0$ 
we define trees which are $(\delta_n)$ 
approximations to $(E_i)$-block trees. 
   
$(G_i)$ is a {\em blocking\/} of $(E_i)$ if there exist integers 
$0 = m_0<m_1<\cdots$ so that $G_i = \oplus_{j=m_{i-1}+1}^{m_i} E_j$ for all $i$. 
$(x_n)\subseteq S_Z$ is a {\em skipped block w.r.t $(G_i)$} if 
\begin{itemize}
\item[(SB)] \qquad for some sequence  $1=k_0<k_1<\cdots < $ in $\N$,
$x_n\in \oplus_{j=k_{n-1}+1}^{k_n-1} G_j$ for all $n$.
\end{itemize}
If $\delta= (\delta_i)$ with $\delta_i\searrow 0$ and $(x_n) \subseteq S_Z$
we say $(x_n)$ is a {\em $(\delta_i)$-skipped block w.r.t.\/} $(G_i)$ if 

\n  ($\delta$-SB)\quad 
for some sequence $1= k_0<k_1<\cdots$ in $\N$, 
$$\|(\Id - P_{\oplus_{j=k_{n-1}+1}^{k_n-1} G_j }) x_n\| 
<\delta_n\ \text{ for all }\ n\ .$$

\begin{thm}\label{T1} 
Let $\cB$ be a countable collection of subsets of $S^{\,\omega}_X$.
Then there exists an isometric embedding of $X$ into a space $Z$ having an
FDD $(E_i)$, so that for $\cA\in\cB$ the following are equivalent.
\begin{enumerate}
\item[a)] $\forall\vp>0 \quad(\text{W}_I(\,\bcAvp\,))$.
\item[b)] For every $\vp>0$ there is a blocking  $(G_i)$ of $(E_i)$
and a sequence $\delta_i\searrow 0$, so that for
every sequence $(x_n)\subset S_X$, satisfying 
$(\delta$-SB) w.r.t. $(G_i)$, $(x_n)\in \bcAvp$.
\item[c)] For every $\vp>0$ there is a blocking  $(G_i)$ of $(E_i)$, so that for
every sequence $(x_n)\subset S_X$ 
$(SB)$ w.r.t. $(G_i)$, $(x_n)\in \bcAvp$.
\end{enumerate}
 
If $X$ has a separable dual $(E_i)$ can be chosen to be shrinking 
and independent from $\cB$, and, furthermore,  if $X$ is reflexive, 
$Z$ can be chosen to be reflexive.  
In these cases (a) is equivalent to 
\begin{enumerate}
\item[d)] 
For every $\vp>0$ every weakly null tree in $S_X$ has a branch in $\bcAvp$.
\end{enumerate}
\end{thm}

\begin{remark} 
Note that Theorem \ref{T1} means the following. 
Assume for all $\vp>0$ Player~I has a winning strategy for the $\bcAvp$-game. 
Then given $\vp>0$,
Player~I can embed $X$ into a space with  an appropriate FDD  
$(F_i)$, and  use the following strategy:
\begin{itemize}
\item[{}]  Take $Y_1=\oplus_{i=2}^\infty F_i\cap X$.
\item[{}]  If Player II has chosen the vector $x_{n-1}$ in the $n-1$st round,
\item[{}]  choose $N\in\N$ so that $\|P_{\oplus_{i=N}^\infty F_i} 
(x_{n-1})\|<\delta_n$ and put 
\item[{}]  $Y_n=\oplus_{i=N+1} F_i\cap X$. 
\end{itemize}
\end{remark}

The proof of Theorem~\ref{T1} also gives the following. 
Suppose $X\subseteq Z$ where $Z$ has an FDD $(E_i)$ and suppose Player~I 
is only allowed to choose subspaces in $\Gamma = \{X\cap \oplus_{i=n}^\infty 
E_i:n \in \N\}$ then a) and b) are equivalent for all $\cA$. 
  
\begin{proof}[Proof of Theorem \ref{T1}]

We first choose a decreasing sequence of finite codimensional spaces 
$(Y_n)$ in $X$  so that for each $\cA\in\cB$ the equivalences 
(a)$\iff$(b)$\iff$(c)$\iff$(d), and, if $X^*$ is separable, (d)$\iff$(e),
of Proposition~\ref{P2} hold. 
Then we choose the space with an FDD $(E_i)$ as in Lemma~\ref{L1}.

We note  that trivially (b) of the statement of Theorem~\ref{T1} implies (c). 
Since the conclusion of Lemma \ref{L1} implies
that every $(X_n)$-block tree (recall, $X_n = \oplus_{i=n+1}^\infty 
E_i\cap X$) has for given sequence 
$\delta_i\searrow 0$ a subtree which  is
a $(\delta_i)$-approximation of an $(E_i)$-block tree for which some 
branch is (SB) w.r.t. $(G_i)$, 
condition~(c) implies condition~(a) (Player~II cannot have a winning 
strategy). 
If $X^*$ is separable the statement (a)$\iff$(d) is exactly the statement of 
the second part of Proposition~\ref{P2}.

Thus, we are left with the verification of the implication (a)$\Rightarrow$(b).

Let $\vp>0$ and $\cA\in\cB$. We put $\eta_i=\vp/c2^{i+2}$, where the 
constant $c>1$ comes from the conclusion of Lemma~\ref{L1} (c).

\begin{claim}
Every tree $(x_A)_{A\in\fN}$ 
in $S_X$ having the property that
\begin{equation}\label{E17}
x_A\in X\cap \bigl(S_{\ds\oplus_{i=\max A+1}^\infty E_i}
\bigr)_{\eta_{\card A}},\text{ whenever }A\in\fN,
\end{equation}
is an $(\vp/2^n)$-approximation to a $(Y_n)$-block tree, and therefore 
must have a branch in $\overline{\cA_{2\vp}}$ 
(Proposition \ref{P2} (a)$\iff$(c)).
\end{claim}
\renewcommand{\qedsymbol}{\ }
\end{proof}

\begin{remark} 
Note that it is in general  not true that if 
$x\in X\cap \bigl(S_{\ds\oplus_{i=m}^n E_i}\bigr)_{\delta}$, then 
we will be able to aproximate $x$ by an element in $X_{m-1}=
 \oplus_{j=m}^\infty E_j\cap X$ up to some 
$r(\delta)$, which converges to $0$
if $\delta$ tends to 0, and which only depends on $\delta$,
but not on $m$ and $n$. But condition~(c) of Lemma~\ref{L1} 
will ensure that we can at least approximate
$x$ by an element of $Y_n$, up to a fixed multiple of $\delta$.
\end{remark}

In order to prove the claim it suffices to show 
\begin{equation}\tag*{$(*)$}
\left\{  \begin{array}{cl}
&\text{Let $\delta >0$ and\ }
x\in X\cap \bigl(S_{\oplus_{i=n+1}^\infty E_i}\bigr)_\delta\ .\\
\noalign{\vskip4pt}
&\text{Then there is a $y\in S_{Y_n}$ with $\|x-y\|\le 4\delta c$.}
\end{array} \right.
\end{equation}

In order to verify the claim we can assume without loss of generality that $\delta<1/2c$ 
(otherwise the claim is trivial).
Choose $u\in\oplus_{i=n}^{\infty} E_i$ and $v\in Z$, $\|v\|<\delta$, 
so that $x=u+v$.
{From} Lemma~\ref{L1}(c) we deduce (recall that 
$D_n\subset \oplus_{i=1}^n E^*_i$) that
$$\|x\|_{X/Y_n}\le c\max_{w^*\in D_n} w^*(x)
= c\max_{w^*\in D_n} w^*(v)< c\delta\ .$$

We can therefore write $x=\tilde y +d$, with $\tilde y\in Y_n$ and 
$d\in X$, satisfying $\|d\|<c\delta$. Since
$\|x\|=1$, we have $1-c\delta<\|\tilde y\|<1+c\delta$. 
Letting $y=\tilde y/\|\tilde y\|$ this implies that $\|x-y\|\le 4c\delta$, 
and finishes the proof of $(*)$.

We next show that there is an increasing sequence $N_i\subset N$ so that
if we let $G_i = \oplus_{s=1+N_{i-1}}^{N_i} E_s$ then  
for every sequence $(x_k)\subset S_X$ for which
there exist integers $m_0=1<m_1<\cdots$ so that
$$
\dist(x_k,\oplus_{s=1+m_{k-1}}^{m_k-1} G_s) =
\dist(x_k,\oplus_{i=1+N_{1+m_{k-1}}}^{N_{m_k-1}} E_i)<\eta_k,\qquad k\in\N,$$
then $(x_k)\in \overline{\cA_{4\vp}}$. 

Since for all $x\in S_X$ it follows that  
($K$ depends on the basis constant of $(E_i)$) 
$$\dist(x,\oplus_{i=m+1}^n E_i)
\le K \|(\Id - P_{\oplus_{i=m+1}^n E_i})(x)\|$$
this will finish the proof of b) taking $\delta_i=\eta_i/K$.

For $\bN=(N_i)_{i=1}^\infty\in\iN$ (the set of infinite subsequences of $\N$) 
we put ($N_0=0$) 
\begin{gather}
F^{\bN}_i=\oplus_{j=1+N_{i-1}}^{N_i} E_i,\qquad i=1,2,\ldots, \text{ and }
\label{E18}\\
\cF^{\bN}=\bigl\{(x_i)_{i=1}^\infty\subset S_X :\forall i\in\N 
\quad\dist(x_i,F^{\bN}_{2i})<\eta_i\bigr\}.
\label{E19}
\end{gather}

\begin{remark} 
For $\bN\in\iN$ and every $(z_i)_{i=1}^\infty\subset S_X$, 
having the property that 
$$\dist(S_{\ds\oplus_{j=1+m_{i-1}}^{m_i-1} F_j^{\bN}}, z_i)<\eta_i,
\qquad i=1,2,\ldots,$$
for some sequence $1\le m_0<m_1<m_1+1<m_2<m_2+1<m_3<\ldots$ there is a sequence
$\bM\in\iNN$ so that $(z_i)\in\cF^{\bM}$.

Indeed, let $\tilde z_i\in S_{\ds\oplus_{j=1+m_{i-1}}^{m_i-1} F_j^{\bN}}$, 
for $i\in\N$ so that $\|\tilde z_i-z_i\|<\eta_i$ and
put $M_{2i-1}=N_{m_{i-1}}$ and $M_{2i}=N_{m_i-1}$. Then it follows that
$$\tilde z_i\in S_{\ds\oplus_{j=1+m_{i-1}}^{m_i-1} F_j^{\bN}}=
S_{\ds\oplus_{s=N_{1+m_{i-1}}}^{N_{m_i-1}} E_s}=S_{F^{\bM}_{2i}}.$$
Thus, $(z_i)\in\cF^{\bM}$.
\end{remark}
      
\begin{proof}[Completion of the proof of Theorem \ref{T1}]

We put 
$$\cC=\{\bN\in\iN: \cF^{\bN}\subset\overline{\cA_{4\vp}}\}\ .$$
It is easy to see that $\cC$ is closed in the pointwise topology on $\iN$, since
$\overline{\cA_{4\vp}}$ is closed with respect to the  product of the 
discrete topology on $S^{\,\omega}_X$.

By the infinite version of Ramsey's theorem (cf.[O]) we deduce that 
one of the following two cases occurs.
\begin{align*}
&\text{Either there exists an }\bN\in\iN\text{ so that }  \iNN\subset\cC\ .\\
&\text{Or there exists an }\bN\in\iN\text{ so that }\iNN\subset\iN
\setminus\cC\ .
\end{align*}
If the first alternative occurs we are finished by the above remark.
Assuming the second alternative, we will show that there is a tree 
in $S_X$ satisfying \eqref{E17} without any
branch in $\overline{\cA_{2\vp}}$. 
This would be a contradiction and imply that 
the second alternative cannot occur.

If we assume the second alternative we can pick for 
each $\bM\in\iNN$ a sequence
$(y^{\bM}_i)_{i=1}^\infty\in\cF^{\bM}$ which is not in $\cC$. 
Let $\bN=\{N_1,N_2,\ldots\}$.

Note that for any $\bM\subset\{N_3,N_4,\ldots\}$, 
$$y_1^{(N_1,N_2,\bM)}\in S_X\cap 
\bigl(S_{\ds\oplus_{i=1+N_1}^{N_2} E_i}\bigr)_{\eta_1}\ .$$
Here $(N_1,N_2,\bM)$ is the infinite sequence starting with $N_1$ 
and $N_2$ and then consisting of the elements of $\bM$).

Using the finite version of Ramsey's theorem and the compactness of 
$S_{\ds\oplus_{i=1+N_1}^{N_2}E_i}$ we can find a vector 
$$x_{\{1\}}\in S_X\cap\bigl(S_{\ds\oplus_{i=1+N_1}^{N_2} E_i}\bigr)_{\eta_1}$$
and an $\bM^{(1)}\subset\{N_3,N_4,\ldots\}$ such that
\begin{equation}\label{E20}
\|x_{\{1\}}-y_1^{(N_1,N_2,\bM)}\|<2\eta_1\text{ for all } 
\bM\in[\bM^{(1)}]^\omega\ .
\end{equation}

Doing the same procedure again, we can find an 
$$x_{\{2\}}\in S_X\cap
\biggl(S_{\ds\oplus_{i=1+N_1^{(2)}}^{N^{(2)}_2} E_i}\biggr)_{\eta_1}$$
and an $\bM^{(2)}\subset[\bM_{(1)}]^\omega$ so that
\begin{equation*}
\|x_{\{2\}}-y_1^{N_1^{(2)},N_2^{(2)},\bM}\|
< 2\eta_1\text{ for all } \bM\in[\bM^{(2)}]^\omega\ ,
\end{equation*}
where $N_1^{(2)}$ and $N_2^{(2)}$ are the first two elements of 
the sequence $\bM^{(1)}$. 
Proceeding this way we construct a sequence
$x_{\{i\}}$ and a decreasing sequence $(\bM^{(i)})$ of infinite 
subsequences of $\bN$ so that
$$x_{\{i\}}\in S_X\cap\bigl(S_{\ds\oplus_{j=1+N_1^{(i)}}^{N^{(i)}_2} E_j}\bigr)_{\eta_1},\text{ and }$$
$$\|x_{\{i\}}-y^{(N^{(i)}_1,N^{(i)}_2,\bM)}\|<\eta_1,
\text{ for all }\bM\in[\bM^{(i)}]^\omega\ .$$
This sequence will be the first level of a tree and the beginning of the 
level by level recursive construction of this tree as follows.

Assume that for some $\ell$ and every $A\in[\N]^{\le\ell}$ we have chosen an
$x_A\in S_X$, a pair of natural number $N_1^{(A)}$, 
and $N_2^{(A)}$, and a sequence
$\bM^{(A)}\in[\{N\in\bN: N>N_2^{(A)}\}]^\omega$ so that the 
following conditions \eqref{E22} and  \eqref{E23} are satisfied.
\begin{align}\label{E22}
&\text{If } A\in[\N]^{<\ell}\cup\{\emptyset\} \text{ and } n>m>\max A \text { then }\\
&\quad N_1^{(A)}\!<\!N_2^{(A)}\!<\!N_1^{(A\cup\{m\})}\!<\!N_2^{(A\cup\{m\})}\!<\!N_1^{(A\cup\{n\})}\!<\!N_1^{(A\cup\{n\})}\notag\\
&\quad[N_1^{(\emptyset)}=N_2^{(\emptyset)}=0]\notag\\
&\quad \bM^{(A)} \supset \bM^{(A\cup\{n\})}\notag
\end{align}
\begin{align}\label{E23}
&\text{If }n_1<n_2<\ldots<n_\ell \text{ are in }\N, \text{ we put }\\
&\qquad\qquad A_j=\{n_1,n_2,\ldots n_j\} \text{ for }j=1,2,\ldots\ell\ .
\notag\\ 
&\text{Then:}\notag\\     
&\qquad x_{A^{(j)}}\in S_X\cap \bigl(S_{\oplus_{s=N_1^{(A_j)}+1}^{N_2^{(A_j)}} 
E_s}\bigr)_{\eta_j}\notag\\
&\qquad\| x_{A^{(j)}}-y^{\displaystyle(N_1^{(A_1)},N_2^{(A_1)},...N_1^{(A_j)},
N_2^{(A_j)},\bM)}\|<\eta_j  \notag\\
& \qquad\qquad\qquad\qquad\text{ whenever }\bM\in[\bM^{(A_j)}]^\omega\notag
\end{align}

Then we can choose for $A\in[\N]^{\ell}$ the elements $x_{A\cup\{1+\max A\}}$, 
$x_{A\cup\{2+\max A\}}$ etc., and the numbers $N_1^{(A\cup\{1+\max A\})}$,
 $N_2^{(A\cup\{1+\max A\})}$,
 $N_1^{(A\cup\{2+\max A\})}$,
 $N_2^{(A\cup\{2+\max A\})}$ etc. and the  sets  $\bM^{(A\cup\{1+\max A\})}$,
 $\bM_1^{(A\cup\{2+\max A\})}$, etc. exactly in the same way we chosed
$x_{\{1\}}$, $x_{\{2\}}$ etc. and the numbers
 $N^{(1)}_1,N^{(1)}_2,N^{(2)}_1,N^{(2)}_2$ etc. for the first level. 

The condition \eqref{E23} implies that for every branch $(z_n)$ of 
the constructed tree there is an $\bM\in[\N]^\omega$ so that 
$\|z_n- y^{\bM}_n|\le 2\eta_n$, for all $n\in\N$.
Since $(y^{\bM}_n)\not\in\overline{\cA_{4\vp}}$ it follows that 
(recall that $\eta_n\le\vp/2^n$) $(z_n)\not\in\overline{\cA_{2\vp}}$, 
which is a contradiction and finishes the proof.
\end{proof}

\section{Subspaces of $(\oplus_{i=1}^\infty F_i)_p$}\label{S3}

The purpose of this section is to use Theorem~\ref{T1} to produce an
intrinsic characterization of a necessary and sufficient 
condition that ensures a given Banach space
$X$ will embed into an $\ell_p$-sum of finite dimensional spaces.

Let $1\le p<\infty$ and let $F_i$  be a finite dimensional space for $i\in\N$.
{\em The $\ell_p$-sum of} $(F_i)$, $(\sum F_i)_p$, 
is the space  of all sequences
$(x_i)$, with $x_i\in F_i$, for $i=1,2\ldots$, so that
$$\|(x_i)\|_p=\biggl(\sum_{i=1}^\infty \|x_i\|_{F_i}\biggr)^{1/p}<\infty\ .$$

\begin{thm}\label{T2}  
Assume that $X$ is reflexive and that there are $1<p<\infty$,
and $C>1$ so that every weakly null tree in $S_X$ has a branch which
is $C$-equivalent to the  unit vector basis of $\ell_p$.

Then $X$ is isomorphic to a subspace of an $\ell_p$-sum of 
finite dimensional spaces.

More precisely, for any $\vp>0$ there exists a finite codimensional  subspace
$\tilde X$ of $X$, so that $\tilde X$ is $(C^2+\vp)$-isomorphic to
a subspace of an $\ell_p$-sum of finite dimensional spaces.
\end{thm}

Before we start the proof, some remarks are in order.

\begin{remark} 
The assumption that $X$ is reflexive is necessary. 
Indeed, James' space  $J$ \cite{Ja1} is not reflexive but  has 
the property that every weakly null
tree in $S_J$ has a branch which is 2-equivalent to the
unit vector basis of $\ell_2$. 
Actually every normalized skipped block with respect
to the shrinking basis of $J$ is 2-isomorphic to the unit vector basis
of $\ell_2$. 
Since every $\ell_2$ sum  of finite dimensional spaces
must be reflexive, $J$ cannot be isomorphic to a subspace of such a space.
\end{remark}
 
In \cite{KW} Kalton and Werner showed a special version of above result. 
They proved  the conclusion of Theorem~\ref{T2} (with $C=1$) 
under the condition that  $X$ does not contain a copy of $\ell_1$ and
every weakly-null type is an $\ell_p$ type.
This means that for every $x\in S_X$ and
every normalized weakly null sequence $(x_n)\subset S_X$ for $t>0$ one has  
\begin{equation}\label{E24}
\lim_{n\to\infty}\|x+tx_n\|=(1+t^p)^{1/p}\ .
\end{equation}
In \cite{KW} it was  shown  that this condition implies
that $X$ must be reflexive, and it is easy to see that it also implies
the hypothesis of Theorem~\ref{T2} with $C=1+\vp$ for any $\vp>0$.

Secondly, let us explain the reason for the $C^2$ term rather than $C$ in 
the conclusion of Theorem~\ref{T2}. 
A normalized basis $(x_i)$ is $C$-equivalent to the unit vector basis of 
$\ell_p$ if there exist constants $A,B$ with $AB \le C$ and 
\begin{equation*}
A^{-1} \biggl( \sum_{i=1}^\infty |a_i|^p\biggr)^{1/p}
\le \Big\| \sum_{i=1}^\infty a_ix_i\Big\| 
\le B\biggl( \sum_{i=1}^\infty |a_i|\biggr)^{1/p}
\tag{$*$}
\end{equation*}
for all scalars $(a_i)$. 
If we had the hypothesis that every weakly null tree in $S_X$ admitted 
a branch $(x_i)$ with this property then we could obtain the conclusion 
of Theorem~\ref{T2} with $C^2$ replaced by $C$. 
However the constants $A,B$ above could vary with each such tree and 
so we can only use $(*)$ with $A$ and $B$ replaced by $C$. 
In this case we only get  $C^2$-embedding into $\ell_p$.


We also note that Kalton \cite{K} proved the following analogous theorem 
for $c_0$: 
Let $X$ be a separable Banach space not containing $\ell_1$. 
If there exists $C<\infty$ so that every weakly null tree in $S_X$ has a 
branch $C$-equivalent to the unit vector basis of $c_0$ then $X$ embeds 
into $c_0$. 

W.B.~Johnson \cite{J2} showed that in the case $X\subseteq L_p$ $(1<p<\infty)$,
if there exists $K<\infty$ so that every normalized sequence in $X$ has a 
subsequence $K$-equivalent to the unit vector basis of $\ell_p$ then $X$ 
embeds into $\ell_p$. 
The tree hypothesis of Theorem~\ref{T2} cannot in general be weakened to the 
subsequence condition as the following example shows. 
(Theorem~\ref{T2} and this example solve some questions raised in 
\cite{J2}.) 

\begin{example}\label{Ex1}
Let $1<p<\infty$. 
There exists a reflexive space $X$ with an unconditional basis so that $X$ 
satisfies:  for all $\vp >0 $ every normalized weakly null sequence in $X$ 
admits a subsequence $1+\vp$-equivalent to the unit vector basis of $\ell_p$. 
Yet $X$ is not
a subspace of an $\ell_p$-sum of finite dimensional spaces.
\end{example}

\begin{proof} 
Fix $1<q<p$. 
We define $X= (\sum X_n)_p$ where each $X_n$ is given as follows. 
$X_n$ will be the completion of $c_{00}([\N]^{\le n})$ under the norm 
$$\|x\|_n = \sup \biggl\{ \biggl( \sum_{i=1}^m \|x|_{\beta_i}\|_q^p 
\biggr)^{1/p} :(\beta_i)_1^m \text{ are disjoint segments in } [\N]^{\le n}  
\biggr\}\ .$$ 

By {\em a segment\/}  we mean a sequence $(A_i)_{i=1}^k\in [\N]^{\le n}$
with  $A_1=\{n_1,n_2,\ldots n_\ell\}$, 
$A_2=\{n_1,n_2,\ldots n_\ell,n_{\ell+1}\}$
$\ldots$ 
$A_k=\{n_1,n_2,\ldots n_\ell,n_{\ell+1}\ldots n_{\ell+k-1}\}$, 
for some 
$n_1<n_2<\ldots n_{\ell+k-1}$. Thus a segement can be seen as an interval of a
branch (with respect to the usual partial order in  $[\N]^{\le n}$), 
while a branch is a maximal segment.

Clearly the {\em node basis\/} $(e_A^{(n)})_{A\in [\N]^{\le n}}$ 
given by $e_A (B) 
= \delta_{(A,B)}$ is a 1-unconditional basis for $X_n$. 
Furthermore the unit vector basis of $\ell_q^n$ is 1-equivalent to 
$(e_{A_i}^{(n)})_1^n$, if $(A_i)_1^n$ is any branch of $[\N]^{\le n} $.

Thus no extension of the tree $(e_A^{(n)})_{A\in[\N]^{\le n}}$
to a weakly null tree of infinite length in $S_X$ has a branch whose
basis distance to the $\ell_p$-unit vector basis is closer than
$\dist_b(\ell_p^{(n)},\ell_q^{(n)})=n^{\frac1q-\frac1p}\to\infty$
for $n\to\infty$.  Since it is clear that in every subspace $Y$
of an $\ell_p$ sum of finite dimensional spaces every weakly
null tree in $S_Y$ must have a branch equivalent (for a fixed
constant) to the  unit vector basis of $\ell_p$ it follows that
$X$ cannot be embedded into 
a subspace of an $\ell_p$-sum of finite dimensional spaces.

Also each $X_n$ is isomorphic to $\ell_p$ and thus $X$ is reflexive. 

It remains to show that if $(x_j)$ is a normalized weakly null  sequence in 
$X$ and $\vp>0$ then a subsequence is $1+\vp$-equivalent to the unit vector 
basis of $\ell_p$. 
By a gliding hump argument it suffices to prove this in a fixed $X_n$. 
We proceed by induction on $n$. 

For $n=1$ the result is clear since $X_1$ is isometric to $\ell_p$. 
Assume the result has been proved for $X_{n-1}$.
By passing to a subsequence and perturbing we may assume that 
$(x_i)_1^\infty$ 
is a normalized block basis of the node basis for $X_n$. 

Let $\vp_i\downarrow 0$ rapidly. 
For $j\in\N$ let $P_j$ be the basis projection of $X_n$ onto 
$[e_A :A\in [\N]^{\le n}$, $\min A=j]$. 
Passing to a subsequence we may assume that $\lim_{i\to\infty} \|P_jx_i\|_n 
= a_i$  and from the definition of $\|\cdot\|_n$ we have $(a_i)_{i=1}^\infty 
\in B_{\ell_p}$. 
Choose $a_0\ge 0$ so that $(a_i)_{i=1}^\infty \in S_{\ell_p}$. 

Passing to a subsequence of $(x_i)$ we may assume that there exist integers 
$1=N_0<N_1<\cdots$ so that 
\begin{enumerate}
\item[(i)] $x_i (\{j\})\ne0 \Rightarrow j\in [N_i,N_{i+1})$ 
\item[(ii)] $P_jx_i=0$ for $j\ge N_{i+1}$ 
\smallskip
\item[(iii)] $\ds  \Big\|\sum_{j\in [N_i,N_{i+1})} P_jx_i\Big\|_n 
= \biggl( \sum_{j\in [N_i,N_{i+1})} \|P_j x_i \|_n^p\biggr)^{1/p}$ 
is within $\vp_i$ of $a_0$. 
\smallskip
\item[(iv)] If $j\in [N_i,N_{i+1})$, $i\ge 1$, then if $a_j\ne 0$,  
$(a_j^{-1}P_jx_\ell)_{\ell>i}$ is $1+\vp_j$-equivalent to the unit 
vector basis of $\ell_p$. 
\item[(v)] If $j\in [N_0,N_1)$ and $a_j\ne0$ then $(a_j^{-1} 
P_jx_\ell)_{\ell=1}^\infty$ is $1+\vp_j$-equivalent to the unit vector 
basis of $\ell_p$.
\smallskip
\item[(vi)] $\ds \biggl(\sum_{N_1}^\infty a_j^p\biggr)^{1/p} <\vp_1$ 
\smallskip
\item[(vii)] If $j\in [N_0,N_1)$ and $a_j =0$ then $\|P_j x_i\|_n\le\vp_i$ 
for all $i$. 
\item[(viii)] If $j\in [N_i,N_{i+1})$ and $a_j =0$ then $\|P_jx_\ell\|_n 
< \vp_\ell$ for $\ell>i$.
\end{enumerate}

Conditions (iv) and (v) use the induction hypothesis and the fact that for 
all $j$,\break
span$(\{e_{\{j\}\cup A)} :A\in [\N]^{n-1}$, $\min A>j\})$ 
is isometric to $X_{n-1}$.  
Our conditions are sufficient to yield (for suitably small $\vp_j$'s) that 
$(x_i)$ is $1+\vp$-equivalent to the unit vector basis of $\ell_p$. 
We omit the standard yet tedious calculations.
\end{proof}

For the proof of Theorem~\ref{T2} we need a result which was shown in 
\cite{KOS}. 
It is based on a trick of W.~B.~Johnson \cite{J2} where part (a) was shown.

\begin{lem}\label{L3} 
{\rm (Lemma 5.1 in \cite{KOS})} 
Let $X$ be a subspace of a space $Z$ having 
a boundedly complete FDD $(F_n)$ and assume
$X$ is w$^*$ closed (since $(F_n)$ is boundedly complete $Z$ is 
naturally a dual space).
Then for all $\vp>0$ and $m\in\N$ there
exists an $n>m$ such that if 
$x= \sum_1^\infty x_i\in B_X$ with $x_i\in F_i$ for all $i$, then 
there exists  $k\in (m,n]$ with 
\begin{enumerate}
\item[a)] $\|x_k\| <\vp$  and 
\item[b)] $\dist (\sum_{i=1}^{k-1} x_i,X)<\vp$.
\end{enumerate}
\end{lem}

\begin{cor}\label{cor:4.4}
Let $X$ be a subspace of the reflexive space $Z$ and let $(F_i)$ be an FDD 
for $Z$. 
Let $\delta_i\downarrow 0$. 
There exists a blocking $(G_i)$ of $(F_i)$ given by 
$G_i = \oplus_{j=N_{i-1} +1}^{N_i} F_j$ for some 
$0= N_0 < N-1 <\cdots$ with the following property. 
For all $x\in S_X$ there exist $(x_i)_1^\infty\subseteq X$ and 
$t_i\in (N_{i-1},N_i]$ for $i\in \N$ so that 
\begin{itemize}
\item[a)] $\ds x=  \sum_{i=1}^\infty x_i$.
\item[b)] For $i\in\N$ either $\|x_i\| <\delta_i$ or 
 $\displaystyle\| P_{\oplus_{j=t_{i-1}+1}^{t_i-1} F_j} (x_i) -x_i\| < \delta_i \|x_i\|$ 
\item[c)] For $i\in\N$, $\|P_{\oplus_{j=t_{i-1}+1}^{t_i-1}F_j} x-x_i\| 
< \delta_i$. 
\end{itemize}
\end{cor}

\begin{proof} 
We choose an appropriate sequence $\vp_i\downarrow0$ depending upon 
$(\delta_i)$ and the basis constant $K$ of $(F_i)$. 
$N_1$ is chosen by the lemma for $\vp = \vp_1$ and $m=1$. 
We choose $N_2 >N_1$ by the lemma for $\vp=\vp_2$ and $m=N_1$ and so on. 

If $x\in S_X$ the lemma yields for $i\in\N$, $t_i\in (N_{i-1},N_i]$ with 
$\|P_{F_{t_i}} (x)\| <\vp_i$ and $z_i\in X$ with 
$\|P_{\oplus_{j=1}^{t_i-1} F_j} (x) - z_i\|<\vp_i$. 
We then let $x_1 =z_1$ and for $i>1$, $x_i = z_i-z_{i-1}$. 
Thus $\sum_{i=1}^n x_i = z_n \to x$ and so a) holds. 

To see c) we note the following 
\begin{equation*}
\|P_{\oplus_{j=t_{i-1}+1}^{t_i-1} F_j} (x) - x_i\| 
 \le \|P_{\oplus_{j=1}^{t_i-1} F_j} (x) - z_i\| 
  + \|P_{\oplus_{j=1}^{t_{i-1}}F_j} (x) - z_{i-1}\| 
< \vp_i + 2\vp_{i-1}\ .
\end{equation*}
Thus 
$$\|P_{\oplus_{j=t_{i-1}+1}^{t_i-1}F_j} (x_i) - x_i\| 
= \|(\Id - P_{\oplus_{j=t_{i-1}+1}^{t_i-1}F_j} ) 
(x_i - P_{\oplus_{j=t_{i-1}+1}^{t_i-1} F_j} x)\| 
< (2K+1) (\vp_i + 2\vp_{i-1})$$
which can be made less than $\delta_i^2$. 
This yields b).
\end{proof} 

\begin{remark} 
The proof yields that the conclusion of the corollary remains valid 
for any further blocking of the $G_i$'s (which would redefine the $N_i$'s). 
\end{remark}

\begin{proof}[Proof of Theorem~\ref{T2}]
We first show that $X$ embeds into $(\sum G_n)_{\ell_p}$ for some 
sequence $(G_n)$ of finite dimensional spaces. 
Then to obtain the $C^2 +\vp$ estimate we adapt an averaging argument 
 similar to the one of \cite{KW}. 

Applying Theorem~\ref{T1} to the set 
$$\cA = \{ (x_i)\in S_X^\omega :(x_i) \text{ is $C$-equivalent to the 
unit vector basis of }\ell_p\}$$
we find a reflexive space $Z$ with an FDD $(F_i)$ with basis constant $K$ 
which isometrically contains $X$ and $\delta_i\downarrow 0$ so that 
whenever $(x_i)\subseteq S_X$ satisfies 
\begin{equation}\label{eq:29}
\|P_{\oplus_{j=n_{i-1}+1}^{n_i-1}F_j} (x_i) - x_i\| < \delta_i
\end{equation} 
for some sequence $1= n_0 < n_1 < \cdots $ in $\N$ it follows that $(x_i)$ 
is $2C$-equivalent to the unit vector basis of $\ell_p$. 
Let $G_i = \oplus_{j=N_{i-1}+1}^{N_i} F_j$ be the blocking given by 
Corollary~\ref{cor:4.4}. 

Let $x\in S_X$, $x= \sum \bar x_i$ with $\bar x_i\in G_i$ for all $i$. 
Choose $(x_i)$ and $(t_i)\subseteq \N$ as in Corollary~\ref{cor:4.4}. 
It follows from \eqref{eq:29} that (for $\delta_i$'s sufficiently small) 
that 
$$(3C)^{-1} \le \Big(\sum \|x_i\|^p\Big)^{1/p} \le 3C$$
and 
$$(4C)^{-1} \le \biggl(\sum_i \|P_{\oplus_{j=t_{i_{-1}}+1}F_j}^{t_i-1} x\|^p
\biggr)^{1/p} \le 4C\ .$$
Let $y_i = P_{\oplus_{j=t_{i-1}+1}^{t_i-1}F_j} x$. 

Since 
$$\frac1{2(K+1)} \max (\|y_i\|,\|y_{i+1}\|) - \delta_i 
\le \|\bar x_i\| \le (2K+1) \|y_i\| +\delta_i$$
it follows that $X$ embeds isomorphically into $(\sum G_i)_{\ell_p} 
\equiv W$. 

We now renorm $W$ so as to contain $X$ isometrically. 
Thus $W$ has $(G_i)$ as an FDD and there exists $\tilde C$ so that 
if $(w_i)$ is any block basis of a permutation of $(G_i)$ then 
\begin{equation}\label{eq:30} 
\tilde C^{-1} (\sum \|w_i\|^p)^{1/p} \le \|\sum w_i\| 
\le \tilde C(\sum \|w_i\|^p)^{1/p}\ .
\end{equation}

We repeat the first part of the proof. 
Let $\vp>0$. 
{From} Theorem~\ref{T1} we may assume that there exist $\delta_i\downarrow0$ 
so that if $(x_i)\subseteq S_X$ satisfies 
\begin{equation}\label{eq:31} 
\|P_{\oplus_{j=n_{i-1}+1}^{n_i-1} G_j} (x_i) - x_i\| <\delta_i
\end{equation} 
for some $1=n_0<n_1<\cdots$ then $(x_i)$ is $C+\vp$-equivalent to the unit 
vector basis of $\ell_p$. 
Moreover we may assume that this is valid for any further blocking of 
$(G_j)$. 
{From} now on we will replace $X$ by the finite codimensional subspace 
$\oplus_{i=2}^\infty G_i\cap X$ and $W$ by $\oplus_{i=2}^\infty G_i$ and replace 
$G_i$ by $G_{i+1}$. 
We will show that this new $X$ can be $C^2+\vp$-embedded into an $\ell_p$ 
sum of finite dimensional spaces. 

Let $H_i = \oplus_{j=N_{i-1}+1}^{N_i} G_j$ be the blocking given 
by Corollary~\ref{cor:4.4}. 
Thus (for appropriately small $\delta_i$'s) from \eqref{eq:31} 
and Corollary~\ref{cor:4.4} we have that  if $x\in S_X$ there exist 
$t_i \in (N_{i-1},N_i]$ so that 
\begin{equation}\label{eq:32}
(C+2\vp)^{-1}\biggl( \sum_{i=1}^\infty \Big\|\sum_{j=t_{m_{i-1}+1}}^{t_{m_i}}
x_j\Big\|^p\biggr)^{1/p} 
\le \|x\| 
\le (C+2\vp) \biggl( \sum_{i=1}^\infty 
\Big\| \sum_{j=t_{m_{i-1}+1}}^{t_{m_i}} x_j\Big\|^p\biggr)
\end{equation}
where $x=\sum x_i$ is the expansion of $X$ w.r.t. the FDD $(G_j)$ for $W$.

Chose $M\in\N$ so that 
\begin{equation}\label{eq:33}
\frac{\tilde C^{2^{1/p}}}M \le \vp\ \text{ and }\  
(C+2\vp)^{-1} - \frac{\tilde C^2}{M^{1/p}} \ge (C+3\vp)^{-1}\ .
\end{equation}

For $i=1,2,\ldots,M$ and $j=0,1,2,\ldots$ set 
$L(i,j) = \oplus_{s= (j-1)M+i+1}^{jM+i-1} H_s \subseteq W$ 
(using $H_n = \{0\}$ if $n\le0$) and let 
$Y_i = (\oplus_{j=0}^\infty L(i,j))_p$. 
Let $Y = (\oplus_{i=1}^M Y_i)_p$. 
We shall prove that $X$ $C^2 +\eta (\vp)$-embeds into $Y$ where $\eta (\vp) 
\downarrow0$ as $\vp\downarrow0$ which will complete the proof.

To do this we first define maps $T_i :X\to Y_i$ for $1\le i\le M$. 
If $x=\sum x_j$ is the expansion of $x$ w.r.t. $(H_j)$ we let 
$$T_i x = \sum_{s=1}^\infty 
\biggl( \sum_{u=(s-1)M+i+1}^{sM+i-1} x_s\biggr) 
\in \biggl( \oplus_{s=1}^\infty L(i,s)\biggr)_p 
= Y_i\ .$$

Let $1\le i\le M$ and $x\in S_X$, $x= \sum x_j$ as above. 
Write $x_j = \sum_{u=N_{j-1}+1}^{N_j} x(j,u)$ as the expansion of $x_j\in H_j$ 
w.r.t. $(G_i)$. 
Let $(t_i) \subseteq \N$ be given by Corollary~\ref{cor:4.4} 
(w.r.t. $(G_j)$). 
{From} several applications of the triangle inequality and 
\eqref{eq:30} and \eqref{eq:32} we have 
\begin{align*}
\|T_i(x)\|&=\left[\sum_{j=0}^\infty \|\sum_{s=(j-1)M+i+1}^{jM+i-1}
x(s)\|^p\right]^{1/p}\\
&\le\left[\sum_{j=0}^\infty
\|\sum_{u=t_{(j-1)M+i}}^{N_{(j-1)N+i}}
x((j-1)M+i,u)+\sum_{s=(j-1)M+i+1}^{jM+i-1} x(s)\right.\\
&\qquad\qquad\qquad\qquad\qquad\left.+\sum_{u=1+N_{jM+i-1}}^{t_{jM+i}}
x(jM+i,u)\|^p\right]^{1/p}\\
&\quad+\!\!\left[\sum_{j=0}^\infty
\|\!\!\sum_{u=t_{(j-1)M+i}}^{N_{(j-1)N+i}} x((j-1)M\!\!+\!\!i,u)
+\!\!\sum_{u=1+N_{jM\!\!+\!\!i-1}}^{t_{jM+i}}
x(jM+i,u)\|^p\right]^{1/p}\\
&\le (C+2\vp)\|x\|+ \left[\sum_{j=0}^\infty \sum_{u=1+N_{jM+i-1}}^{N_{jM+i}}
\|x(jM+i,u)\|^p \right]^{1/p}\\
&\le (C+2\vp)\|x\|+\tC\|\sum_{s=i(\mod M)} x_s\|\ .
\end{align*}

Similarly one has 
$$\| T_ix\| \ge (C+2\vp)^{-1} \|x\| - \tC \Big\|\sum_{s=i(\mod M)} x_s\|\ .$$

Finally we define $T:X\to Y = (\sum_1^M Y_i)$, by 
$Tx = \frac1{M^{1/p}} \sum_{i=1}^M T_ix$. 
Note that 
\begin{align*}
\|Tx\| & \le \frac1{M^{1/p}} (C+2\vp) 
\biggl( \sum_{i=1}^M \|x\|^p\biggr)^{1/p} 
+ \frac{\tC}{M^{1/p}} 
\biggl( \sum_{i=1}^M \Big\|\sum_{j=i(\mod M)} x_j\Big\|\biggr)\\
&\le (C+2\vp) \|x\| + \frac{\tC^2}{M^{1/p}} \|x\| 
< (C+3\vp) \|x\|
\end{align*}
using \eqref{eq:30} and \eqref{eq:33}. 

Similarly one deduces that for $x\in X$ it follows that 
 $\|T(x)\|\ge \frac1{C+3\vp}\|x\|$.
\end{proof}

\begin{remark} 
The proof of Theorem \ref{T2} had two steps. 
In the first we started with an embedding of $X$ into a certain 
reflexive space $Z$ with an FDD $(F_i)$ and showed that $(F_i)$ can be blocked
to an FDD $(G_i)$ so that $X$ is isomorphic to a subspace of 
$(\oplus G_i)_{\ell_p}$. In that step we could not
deduce any bound for the constant of that isomorphism. 
In the second step we  ``inflated'' $(\oplus G_i)_{\ell_p}$ to the space
$(\oplus_{i=1}^M\oplus_{j\not=i(\mod M)} G_j)_{\ell_p}$ and showed that 
this space contains
a finite codimensional subspace which is $C^2+\vp$-equivalent to $X$.
   
The following example  shows that even if the space 
$X$ has a basis to begin with, it is in general not possible
to pass to a blocking $(F_n)$ of that basis and deduce that 
for some $n_0$ the identity is a $C^2+\vp$-isomorphism between
$\oplus_{n=n_0}^\infty F_n$ and $(\oplus_{n=n_0}^\infty F_n)_{\ell_p}$. 
\end{remark}
    
\begin{example} 
Let $\cD$ be the set of all  sequences $(D_n)$ of pairwise 
disjoint subsets of $\N$, 
so that  for each $n\in\N$, $D_n$ is either a singleton
or it is of the form $D_n=\{k,k+1\}$ for some $k\in\N$.
We give $\ell_2$ the following equivalent norm $|||\cdot|||\ :$
$$|||x|||=\sup_{(D_n)\in\cD} \left(\sum_{n=1}^\infty
\Biggl(\sum_{j\in D_n} |x_j|\Biggr)^2\right)^{1/2},$$
whenever $x=(x_j)\in\ell_2$.
     
It is easy to see that every  normalized skipped block 
$(x^{(n)})$ in $X=(\ell_2,|||\cdot|||)$ is
isometrically eqivalent to the $\ell_2$ unit vector basis. 
Thus the assumptions of Theorem \ref{T2} are satisfied for any $C>1$. 
On the other hand for any blocking $(F_n)$ of the unit vector basis 
$(e_i)$ of $X$ it follows
for any $n$  and $N_n=\max\{N|e_N\in F_n\}$ that 
$e_{N_n+1}\in F_{n+1}$ and that 
the span of $e_{N_n}$ and $e_{N_n+1}$ is isometric to $\ell^2_1$. 
Therefore the norm of the identity  between 
$(\oplus_{n=2}^\infty F_n)_{\ell_2}$ and 
$(\oplus_{n=2}^\infty F_n)_{\ell_2,|||\cdot|||}$ is at least $\sqrt 2$. 
\end{example}
    
The following result shows that the property that every normalized
weakly null tree contains a branch which is $C$-equivalent to 
the $\ell_p$ unitvector basis dualizes.
It can be seen as the isomorphic version of Theorem 2.6. in \cite{KW}. 
     
\begin{cor}\label{cor:4.5}
Assume $X$ is a reflexive Banach space.
For $1<p<\infty$ and $\frac1p +\frac1q =1$ the following
statements are equivalent.
\begin{enumerate}
\item[a)] There is a $C\ge 1$ so that every normalized weakly null
tree in $X$ has a branch which is $C$-equivalent to the  unit
vector basis of $\ell_p$.
\item[b)] There is a $C\ge 1$, a finite codimensional subspace $\tX$ of $X$, 
a sequence of finite dimensional spaces $(E_i)_{i=1}^\infty$,
and an operator $T:\tX\to(\oplus_{i=1}^\infty E_i)_{\ell_p}$, so that
$C^{-1}\|x\|\le \|T(x)\|\le C\|x\|$ for all $x\in \tilde X$.
\item[c)] There is a $C\ge 1$ so that every normalized weakly null
tree in $X^*$ has a branch which is $C$-equivalent to the  unit
vector basis of $\ell_q$.
\item[d)] There is a $C\ge 1$, a finite codimensional subspace $Y$ of $X^*$, 
a sequence of finite dimensional spaces $(E_i)_{i=1}^\infty$,
and an operator $T:Y\to(\oplus_{i=1}^\infty E_i)_{\ell_q}$, so that
$C^{-1}\|x\|\le \|T(x)\|\le C\|x\|$ for all $x\in Y$. 
\end{enumerate}
\end{cor}

\begin{proof}
The implications  (a)$\Rightarrow$(b) and (c)$\Rightarrow$(d) follow
from Theorem~\ref{T2} and its proof. 
If we prove 
(b)$\Rightarrow$(c) then (d)$\Rightarrow$(a) will follow. 
    
Assume that $C\ge 1$, $(E_i)_{i=1}^\infty$, $\tX\subset X$ and
$T:\tX\to Z=(\oplus_{i=1}^\infty  E_i)_{\ell_p}$ are given as in the statement
of (b). By passing to the renorming $|||\cdot|||$,
$|||x|||=\|T(x)\|$, for $x\in\tX$ we can assume without
loss of generality that $\tX$ is isometric to a subspace of $Z$. 
    
We will show that $\tX^*$ satisfies the condition (c). Since 
$\tX^*$ is isomorphic to a subspace of $X^*$ of finite codimension
the claim will follow.
    
Thus let $E:\tX\to (\oplus_{i=1}^\infty E_i)_{\ell_p}$  be an isometric
embedding and let 
$(x^*_A)_{A\in\fN}$ be a normalized weakly null tree in $\tilde X^*$.
    
We will need the following observation.
    
\begin{claim} 
If $(x_n^*)$ is a normalized and weakly null sequence in $\tX^*$,
then there are normalized weakly null sequences $(z^*_n)$ and $(x_n)$ in
$Z^*$ and  $\tX$ respectively so that,
$E^*(z^*_n)=x^*_n$ and $x^*_n(x_n)=1$ for $n\in\N$.
\end{claim}
    
To see this use the Hahn-Banach theorem to     
choose a  normalized sequence $(z^*_n)_{n\in\N}$ in $Z^*$ so that
$E^*(z_n^*)=x_n^*$. 
The sequence  $(z^*_n)$ is weakly null. 
Indeed, otherwise we could
choose a $y^*\in Z^*$, $y^*\not=0$, a subsequence $(z^*_{n_k})$ and
a  weakly null sequence $(y^*_k)$ in $Z^*$ so that
$z^*_{n_k}=y^*+y^*_k$ for all $k\in\N$. 
Thus, $x^*_{n_k}=E^*(y^*)+E^*(y^*_k)$, which implies that $E^*(y^*)=0$ and
therefore that $E^*(y^*_k)=x^*_{n_k}$.
Since $\limsup_{k\to\infty}\|y^*_k\|=
\limsup_{k\to\infty}(\|z^*_{n_k}\|^q-\|y^*\|^q)^{(1/q)}<1$, we get 
a contradiction.
      
Then we choose  $(x_n)\in \tX$ so that $x^*_n(x_n)=1$.  
By a similar argument we have that $(x_n)$ is also weakly null.
     
Using the claim we can find a normalized weakly null tree 
$(z^*_A)_{A\in\fN}$ in
$Z^*$ and a normalized weakly null tree $(x_A)_{A\in\fN}$ in $\tX$,
so that  $E^*(z^*_A)=x^*_A$ and $x^*_A(x_A)=1$ for $A\in\fN$.
      
Given an $\vp>0$ we can choose a branch $(x^*_n)=(x^*_{A_n})$ so that
$(z^*_{A_n})$ is $(1+\vp)$ equivalent  to the  unit vector basis of $\ell_q$,
and $(x_{A_n})$ is $(1+\vp)$ equivalent  to the unit vector basis of $\ell_p$.
This easily implies that $(x^*_{A_n})$ is $(1+\vp)$ equivalent to the
unit vector basis of $\ell_q$.     
\end{proof} 

\begin{remark}
W.B. Johnson and M. Zippin \cite{JZ} proved the following. 
Let $C_p = (\oplus_{i=1}^\infty E_i)_{\ell_p}$ where $(E_i)$ is dense, in the 
Banach-Mazur sense, in the set of all finite dimensional spaces. 
Then $X$ embeds into $C_p$ if and only if $X^*$ embeds into $C_q$ 
(where $\frac1p +\frac1q =1$). 
Thus Corollary~\ref{cor:4.5} could be deduced from \cite{JZ} and 
Theorem~\ref{T2} (and \cite{JZ} could be deduced from the corollary and 
theorem). 

Furthermore the proof of Corollary~\ref{cor:4.5} yields some quantitative 
information. 
If a) holds then b) is true with $C$ replaced by $C+\vp$ for any $\vp>0$. 
If b) holds then c) is valid with $C$ replaced by $C^2+\vp$.  
\end{remark}
    
\section{Spectra and asymptotic structures }\label{S4}
 
In \cite{Mi}  Milman introduced the notion of the 
{\em spectra of a function\/} defined on $S^n_X$. 
Let $(M,\rho)$ be a compact metric space and let $f:S^n_X\to M$ be Lipschitz.
$\sigma(f)$ is defined to be the set of all $a\in M$ for which the 
following condition \eqref{E39} is true
\begin{align}\label{E39}
\forall\, \vp>0\, \forall\,  Y_1\in\cof(X)\, \exists\, 
& y_1\in S_X\,\forall\, Y_2\in\cof(X)\,\exists\, y_2\in S_X\\
&\ldots\,\forall\, Y_n\in\cof(X)\,\exists\, y_n\in S_X \text{ so that }\notag\\
&\rho(f(y_1,y_2,\ldots y_n),a)<\vp \notag
\end{align}
  
In terms of the game we introduced in  Section \ref{S1},  
$\sigma(f)$ is the set of  all  $a\in M$ so that
for any $\vp>0$ Player II has a winning strategy in the $\cA^\vp$-game, where
$$\cA^\vp=\{(y_i)_{i=1}^n\in S^n_X : \rho(a,f(y_1,\ldots y_n))>\vp\}$$
(which means that Player II is able to get 
$f(y_1,\ldots y_n)$ arbitrarily close to $a$).
  
As mentioned in \cite{Mi} one can also  define the spectrum relative 
to any {\em filtration\/} $\cS\subset\cof(X)$, 
meaning that $\cS$ has 
the property that if $X,Y\in\cS$ there is a $Z\in\cS$ for 
which $X\cap Y\supset Z$.
The spectrum of $f$ relative to $S$ is the set $\sigma(f,\cS)$ 
of all $a\in M$ for which
\begin{align}\label{E40}
&\,\forall\,\vp>0\,\forall\, Y_1\in\cS\,\exists\, y_1\in S_{Y_1}
\,\forall\, Y_2\in\cS\,\exists\, y_2\in S_{Y_2}\ldots\,\forall\, 
Y_n\in\cS\,\exists\, y_n\in S_{Y_n}\text{ with}\\
&\rho(f(y_1,y_2,\ldots y_n),a)<\vp. \notag
 \end{align}

It is obvious that $\sigma(f,\cS)\subset \sigma(f,\tilde{\cS})$ 
whenever $\tilde{\cS}\subset\cS$. In particular
it follows that $\sigma(f)\subset\sigma(f,\cS)$ for any filtration $\cS$.
  
If $X$ is a subspace of a space $Z$ with  FDD $(E_i)$ 
we can consider the filtration
$\cS=\{X\cap\oplus_{i=n}^\infty F_i:n\in N\}$ and 
we write $\sigma(f,(F_i))=\sigma(f,\cS)$.
  
On one hand the unrelativized spectrum $\sigma(f)$ seems to be 
the right concept to study
geometric  and structural properties of $X$, since it is ``coordinate free''. 
On the other hand spectra with respect to an FDD is 
combinatorically easier to use and understand.  

But from Theorem \ref{T1} we deduce that $\sigma(f)$ is equal to the spectrum 
with respect to a  certain FDD (of some super space).
 
\begin{prop}\label{P3}
Let $f: S^n_X\to M$ be Lipschitz. 
Then 
\begin{equation}\label{E41}
\sigma(f)=\bigcap\{ C: C\text{ is a closed subset 
of $M$ and } (\text{W}_{I}(f^{-1}(C))) \}\ .
\end{equation}
Moreover for any $\vp>0$, $(\text{W}_{I}(f^{-1}(C))_\vp)$. 
 
Furthermore $X$ can be embedded into a space $Z$ with FDD $(F_i)$ so that for 
every $\vp>0$ there is a $\delta>0$ and an $M_0\in\N$ with the following
property.
 
Whenever  $M_0<M_1<M_2<\ldots M_n$ and  
$(x_i)_{i=1}^n \subseteq S_X$ satisfies 
$$d\Big(x_i, S_{\oplus_{j=1+M_{i-1}}^{M_i-1} F_j}\cap X \Big)
<\delta \text{ for }i=1,\ldots ,n$$
then $\rho(f(x_1,x_2,\ldots x_n),\sigma(f))<\vp$.
 
In the case that $X^*$ is separable, $\sigma(f)$ is the minimal closed 
subset of $M$ so that for any $\vp>0$ any weakly null tree in $S_X$ 
of length $n$ has a branch $(x_1,\ldots x_n)$
so that $\rho(f(x_1,\ldots x_n),\sigma (f))<\vp$.
\end{prop}
 
\begin{proof}
Let $\cC$ denote the set of all closed subsets of $M$ for
which  (W$_{I}(f^{-1}(C))$ holds. 
For $a\in M$ we denote the $\vp$-neighborhood by $U_\vp(a)$
and observe the following equivalences
\begin{align*}
a&\not\in\sigma(f)\\
&\iff\exists\, \vp>0\,\exists\, Y_1\in\cof(X)\,\forall\, y_1\in S_{Y_1}\ldots
\,\exists Y_n\in\cof(X)\,\forall\, y_n\in S_{Y_n}\\
  &\qquad\qquad \rho(f(y_1,..y_n),a)>\vp\\
  &\iff \exists\, \vp>0 \quad (\text{W}_{I}(f^{-1}(M\setminus U_\vp(a)))  \\
  &\iff \exists\, C\in\cC\,,\qquad a\not\in C.
  \end{align*}
 
Thus $\sigma(f)=\bigcap \{C:C\in \cC\}$. 
If $\eta>0$ then $M\setminus (\sigma(f))_\eta$ is compact and is 
contained in the open covering $\bigcup_{C\in\cC} M\setminus C$.
Thus there exists a finite $\tilde{\cC}\subset\cC$ so that  
$M\setminus (\sigma(f))_\eta\subset \bigcup_{C\in\tilde{\cC}} M\setminus C$ 
and thus
$(\sigma(f))_\eta\supset\bigcap_{C\in\tilde{\cC}} C$ 
which implies by Proposition \ref{P0} that Player I has a winning strategy 
for $f^{-1}((\sigma(f))_\eta)$. 
By the uniform continuity of $f$,  
$\eta$  can be chosen small enough so that
$f^{-1}((\sigma(f))_\eta)_n$ contained in a given neighborhood of 
$f^{-1}(\sigma(f))$ which finishes the proof of the first part.
The remainder of the proposition follows easily from Theorem~\ref{T1}. 
\end{proof}

A special example of spectra was considered by Milman and Tomczak 
\cite{MT}, the {\em asymptotic structure of $X$}. 
A finite dimensional space $E$  together with a normalized monoton basis
$(e_i)_1^n$ is called an {\em element of the $n^{th}$-asymptotic 
structure of $X$}  and we write $(E,(e_i)_{i=1}^n)\in\{X\}_n$ if
\begin{align}\label{E42}
&\forall\,\vp>0\,\forall\, Y_1\in\cof(X)\,\exists\, y_1\in S_X\ldots 
\,\exists\, Y_n\in\cof(X)\,\exists\, y_n\in S_X\\
&\dist_b((y_i)_{i=1}^n,(e_i)_{i=1}^n)<1+\vp \notag
\end{align}
where $\dist_b$ denotes the basis distance, i.e.,
if $(e_i)_{i=1}^n$ and $(f_i)_{i=1}^n$ are
two  bases of $E$ and $F$ respectively then
$\dist_b((e_i)_{i=1}^n,(f_i)_{i=1}^n)$ is defined to be 
$\|T\|\cdot\|T^{-1}\|$ where $T:E\to F$ is given by 
$T(e_i)=f_i$, for $i=1,\ldots n$. 
Note that the space $(M_n,\log \dist_b)$ of all normalized 
bases of length $n$ and basis constant not exceeding a fixed 
constant is a compact metric space.
 
Therefore we deduce from Proposition~\ref{P3} and 
the usual diagonalization argument the following Corollary 
(cf. \cite{KOS}). 

\begin{cor}\label{C1}
$X$ can be embedded into a space $Z$ with FDD $(F_i)$ so that for 
every $k\in\N$ it follows that:
 
Whenever  $k=M_0<M_1<M_2<\ldots M_k$ and 
$$x_i\in S_{\oplus_{j=1+M_{i-1}}^{M_i-1} F_j}\cap X 
\text{ for }i=1,2\ldots k\ ,$$
then $\dist_b((x_i)_{i=1}^k,\{X\}_k)<1+\vp$.
 
In the case that $X^*$ is separable, $\{X\}_k$ is 
the minimal closed subset of $M_k$ so 
that for any $\vp>0$ any weakly null tree in $S_X$ of 
length $n$ has a branch $(x_1,\ldots x_k)$
so that $\dist_b((x_i)_{i=1}^k,\{X\}_k)<1+\vp$.
\end{cor}

An interesting case is when the asymptotic structure of $X$ is as 
small as possible.

\begin{thm}\label{XsepreflexBS}
Let $X$ be a separable reflexive Banach space with $|\{X\}_2|=1$. 
Then there exists $p\in (1,\infty)$ so that $X$ embeds into the $\ell_p$-sum 
of finite dimensional spaces. 
Moreover for all $\vp>0$ there exists a finite codimensional subspace 
$X_0$ of $X$ which $1+\vp$-embeds into the $\ell_p$-sum of finite 
dimensional spaces. 
\end{thm}

\begin{proof} 
Since there exists $1\le p\le\infty$ so that the unit vector basis of 
$\ell_p^2$ is in $\{X\}_2$ (see \cite{MMT}) we have that $\{X\}_2$ must 
be this unit vector basis. 
In turn this condition (see \cite{MMT} or \cite{KOS}) implies that $X$ 
contains an isomorph of $\ell_p$ ($c_0$ if $p=\infty$) and so $1<p<\infty$. 

Let $X\subseteq Z$, a reflexive space with an FDD $(E_n)$. 
The condition on $\{X\}_2$ yields that for all $\vp>0$ there exists $n$ 
so that if $x_1 \in S_X\cap [E_i]_{i=n}^\infty$ then there exists $m$ 
so that if $x_2\in S_X \cap [E_i]_{i=m}^\infty$ then $(x_i)_1^2$ is 
$1+\vp$-equivalent to the unit vector basis of $\ell_p^2$. 
{From} this it follows that $X$ satisfies the hypothesis of 
Theorem~\ref{T2} with $C=1$ and thus the theorem follows.
\end{proof}

The following problem remains open. 
We say $X$ is {\em Asymptotic\/} $\ell_p$ if there exists $K<\infty$ so that 
for all $k$ and all $(x_i)_1^k \in \{X\}_k$, $(x_i)_1^k$ is $K$-equivalent 
to the unit vector basis of $\ell_p$. 
An FDD $(E_n)$ for a space $Z$ is {\em asymptotic\/} $\ell_p$ if there 
exists $K<\infty$ so that for all $k$ if $(x_i)_1^k$ is a block sequence 
of $(E_i)_k^\infty$ in $S_Z$, then $(x_i)_1^k$ is $K$-equivalent to the 
unit vector basis of $\ell_p$. 

\begin{prob}\label{prob:5.4}
Let $X$ be a reflexive Asymptotic $\ell_p$ space for some $1<p<\infty$. 
Does $X$ embed into a space $Z$ with an asymptotic $\ell_p$ FDD?
\end{prob}


\begin{thebibliography}{MMT}

\bibitem[Ja1]{Ja1}
R.C. James, 
{\em Uniformly nonsquare Banach spaces},
Ann. of Math. (2){ \bf 80} (1964), 542--550.

\bibitem[J2]{J2}
W.B. Johnson,
{\em On quotients of $L_p$ which are quotients of $\ell_p$},  
Compositio Math. {\bf 34} (1977),  69--89. 

\bibitem[JZ]{JZ}
W.B. Johnson and M. Zippin, 
{\em Subspaces and quotient spaces of $(\sum G_n)_{\ell_p}$ and 
$(\sum G_n)_{c_0}$}, 
Israel J. Math. {\bf 17} (1974), 50--55. 

\bibitem[K]{K}
N.J. Kalton,  
{\em On subspaces of c$_0$ and extensions of operators into
C$(K)$-spaces}, preprint

\bibitem[KW]{KW}
N.J. Kalton and D. Werner, 
{\em Property $(M)$, $M$-ideals, and almost isometric structure 
of Banach spaces}, 
J. Reine und Angew. Math. {\bf 461} (1995), 137--178.

\bibitem[KOS]{KOS}
H.~Knaust, E.~Odell, and Th.~Schlumprecht, 
{\em On asymptotic structure, the Szlenk index and UKK properties in 
Banach spaces},
Positivity {\bf3} (1999), 173--199.

\bibitem[Ma]{Ma} 
D.A. Martin,  
{\em Borel determinacy}, 
Annals of Math. {\bf 102}  (1975),  363--371. 
 
\bibitem[MMT]{MMT}
B. Maurey, V.D. Milman and N. Tomczak-Jaegermann, 
{\em Asymptotic infinite-dimensional theory of Banach spaces}, 
Oper. Theory: Adv. Appl. {\bf 77} (1994), 149--175. 
 
\bibitem[Mi]{Mi} 
V. Milman, 
{\em Geometric theory of Banach spaces II, geometry of the unit sphere}, 
Russian Math. Survey {\bf 26} (1971), 79--163 
(translation from Russian).

 
\bibitem[MT]{MT}
V.D. Milman and N. Tomczak-Jaegermann, 
{\em Asymptotic $\ell_p$ spaces and bounded distortions},
eds. Bor-Luh Lin and W.B. Johnson,
Contemp. Math. {\bf 144} (1993), 173--195.

\bibitem[O]{O} 
E. Odell, 
{\em Applications of Ramsey theorems to Banach space theory}, 
Notes in Banach spaces, 
ed. H.E. Lacey, 
Univ. of Texas Press, 
Austin, TX 
(1980), 
379--404 

\bibitem[W]{W} 
R. Wagner,
{\em Finite high-order games and an inductive approach towards Gowers' 
dichotomy}, 
Annals of Pure and Applied Logic, to appear. 

\bibitem[Z]{Z} 
M. Zippin, 
{\em Banach spaces with separable duals}, 
Trans. AMS, {\bf 310}, Nr. 1 (1988), 371--379.




\end{thebibliography}
 \end{document}